\documentclass[12pt,reqno]{amsart}
\usepackage{amsthm,amsfonts,amssymb,amscd,amsmath,euscript}
\usepackage[all]{xy}
\newtheorem{theorem}[subsection]{Theorem}
\newtheorem*{thm}{Theorem}

\newtheorem{proposition}[subsection]{Proposition}

\newtheorem{lemma}[subsection]{Lemma}
\newtheorem{corollary}[subsection]{Corollary}

\theoremstyle{definition}
\newtheorem{definition}[subsection]{Definition}

\newtheorem{proposition-definition}[subsection]{Proposition-Definition}

\theoremstyle{remark}
\newtheorem{remark}[subsection]{Remark}

\newcommand{\actson}{\begin{picture}(5,4)\thinlines
\put(2,1.6){\oval(4,2.5)[r]}
\put(2,2.8){\vector(-1,0){0.5}}\end{picture}
}

\newcommand{\ccd}{{\begin{picture}(2,2)
\put(1,1){\circle*{0.7}}\end{picture}}}
\newcommand{\gseven}{\Gamma_{12}^{7}}

\newcommand{\dual}{{\scriptscriptstyle \vee}}
\newcommand{\NA}{{\overline{NA}}}
\newcommand{\NE}{{\overline{NE}}}
\newcommand{\sten}{{\Sigma_{12}^{10}}}
\newcommand{\XT}{{X_{12}}}

\newcommand{\Ext}{\operatorname{Ext}\nolimits}

\newcommand\Hilb{{\operatorname{Hilb}\nolimits}}

\newcommand{\Hom}{\operatorname{Hom}\nolimits}
\newcommand{\id}{\operatorname{id}\nolimits}
\newcommand{\length}{\operatorname{length}\nolimits}

\newcommand{\mult}{\operatorname{mult}\nolimits}

\newcommand{\Pic}{\operatorname{Pic}\nolimits}
\newcommand{\pr}{\operatorname{pr}\nolimits}

\newcommand{\rk}{\operatorname{rk}\nolimits}

\newcommand{\Sing}{\operatorname{Sing}}

\newcommand{\Spin}{\operatorname{Spin}(10)}
\newcommand{\Supp}{\operatorname{Supp}}

\newcommand{\CC}{{\mathbb C}}
\newcommand{\RR}{{\mathbb R}}
\newcommand{\ZZ}{{\mathbb Z}}
\newcommand{\QQ}{{\mathbb Q}}
\newcommand{\PP}{{\mathbb P}}
\newcommand{\NN}{{\mathbb N}}
\newcommand{\FF}{{\mathbb F}}

\newcommand{\BBB}{{\mathcal B}}
\newcommand{\OOO}{{\mathcal O}}

\newcommand{\III}{{\mathcal I}}
\newcommand{\JJJ}{{\mathcal J}}

\newcommand{\EEE}{{\mathcal E}}

\newcommand{\LLL}{{\mathcal L}}
\newcommand{\FFF}{{\mathcal F}}
\newcommand{\CCC}{{\mathcal C}}

\newcommand{\NNN}{{\mathcal N}}
\newcommand{\MMM}{{\mathcal M}}

\newcommand{\UUU}{{\mathcal U}}

\newcommand{\EXT}{{\mathcal Ext}}

\newcommand{\cal}[1]{\mathcal{#1}}
\renewcommand{\bar}[1]{\overline{#1}}

\newcommand\alp{\alpha}

\newcommand\ep{\varepsilon}

\renewcommand\phi{\varphi}

\newcommand{\isoto}{{\lra\hspace{-1.3 em}
\raisebox{ 0.6 ex}{$\textstyle\sim$}\hspace{0.8 em}}}
\newlength{\rrrr}
\newcommand{\into}{\hookrightarrow}

\newcommand\lra{{\longrightarrow}}

\newcommand\rar{\rightarrow}
\newcommand\lrdash{\:
\xymatrix@1{\ar@{-->}[r]&}\:
}
\newcommand{\lrdashar}[1]{\:
\xymatrix@1{\ar@{-->}[r]^{#1}&}\:
}

\newcommand{\simvert}{\mbox{\raisebox{0.8 ex}
{{\scriptsize )}\hspace{-1.ex}\raisebox{.7em}{\scriptsize (}}\rule{0pt}{1pt}}}

\renewcommand\square{\frame{\phantom{{\large x}}}}
\renewcommand\emptyset{\varnothing}
\newcommand\empt{\varnothing}

\author{A. Iliev$^1$}

\address{\scriptsize Atanas Iliev: 
Institute of Mathematics,
Bulgarian Academy of Sciences,
Acad. G. Bonchev Str., 8,\ 
1113 Sofia, Bulgaria}
\email{ailiev@math.bas.bg}

\author{D. Markushevich$^2$}
\address{\scriptsize Dimitri Markushevich:  
Math\'ematiques - b\^{a}t.M2, Universit\'e Lille 1, 
F-59655 Villeneuve d'Ascq Cedex, France}
\email{markushe@agat.univ-lille1.fr}


\title{Parametrization of Sing$\,\boldsymbol{\Theta}$ for a Fano 3-fold of genus 7
by moduli of vector bundles}

\begin{document}

\begin{abstract}
According to Mukai, any prime
Fano threefold $X$ of genus 7 is a linear section of the
spinor tenfold in the projectivized half-spinor space
of Spin(10). The orthogonal linear section of the spinor
tenfold is a canonical genus-7 curve $\Gamma$, and
the intermediate Jacobian $J(X)$ is isomorphic to the Jacobian
of $\Gamma$. It is proven that, for a generic $X$, the Abel-Jacobi map of
the family of elliptic sextics on $X$ factors through
the moduli space of rank-2 vector
bundles with $c_1=-K_X$ and $\deg c_2=6$ and that the
latter is birational
to the singular locus of the theta divisor of $J(X)$.
\end{abstract}

\maketitle

\footnotetext[1]{Partially supported by the grant MI-1503/2005
of the Bulgarian Foundation for Scientific Research}

\footnotetext[2]{Partially supported by the grant
INTAS-OPEN-2000-269}

\setcounter{section}{-1}

\section{Introduction}

This work is a sequel to the series of papers on moduli spaces
$M_X(2;k,n)$ of stable rank-2 vector
bundles on Fano 3-folds $X$ with Picard group $\ZZ$
for small Chern classes $c_1=k,\ c_2=n$.
The nature of the results depends strongly
on the index of $X$, which is defined as
the largest integer that divides the canonical class $K_X$
in $\Pic X$. Historically, the first Fano 3-fold for which
the geometry of such moduli spaces was studied
was the projective space
$\PP^3$, the unique Fano 3-fold of index 4. 
The most part of results for $\PP^3$
concerns the problems of rationality,
irreducibility or smoothness of the moduli space, see
\cite{Ba-1}, \cite{Ba-2}, \cite{Ha}, \cite{HS}, \cite{LP}, 
\cite{ES}, \cite{HN}, \cite{M}, \cite{BanM},
\cite{GS}, \cite{K}, \cite{KO}, \cite{CTT}
and references therein.

The next case is
the $3$-dimensional quadric $Q^3$, which is Fano of
index~3. Much less is known here, see \cite{OS}.
Further, the authors
of \cite{SW} identified the moduli spaces $M_X(2;-1,2)$
on all the Fano 3-folds $X$ of index 2 except for the
double Veronese cone $V'_1$,
which are (in the notation of Iskovskikh)
the quartic double solid $V_2$, a 3-dimensional cubic $V_3$,
a complete intersection of two quadrics $V_4$,
and a smooth 3-dimensional section of the Grassmannian
$G(2,5)$ by three hyperplanes $V_5$.
It turns out that all the vector bundles in
$M_X(2;-1,2)$ for these threefolds
are obtained by Serre's construction from
conics. Remark that for $\PP^3$ and $Q^3$ all the known
moduli spaces are either rational or supposed to be rational,
whilst \cite{SW} provides first nonrational examples.

We will also mention the paper \cite{KT} on the moduli of
stable vector bundles on
the flag variety $\FF (1,2)$, though it is somewhat apart, for
$\FF (1,2)$ has Picard group $2\ZZ$. This is practically all
what was known on the subject until the year 2000, when
a new tool was brought into the study of the moduli spaces:
the Abel--Jacobi map
to the intermediate Jacobian $J(X)$. For the 3-dimensional
cubic $X=V_3$, it was proved in \cite{MT-1}, \cite{IM-1} that
the open part of 
$M_X(2;0,2)$ parametrizing the vector bundles obtained by Serre's
construction from elliptic quintics is sent by the Abel--Jacobi
map isomorphically onto an open subset of $J(X)$. Druel
\cite{D} proved the irreducibility of $M_X(2;0,2)$ and
described its compactification by semistable sheaves;
see also the survey \cite{Beau-1}. The other index-2 case,
that of the double solid $V_2$, was considered in \cite{Ti},
\cite{MT-2}, where it was proved that the vector bundles
coming from the elliptic quintics on $V_2$ form an irreducible
component of $M_{V_2}(2;0,3)$ on which the Abel--Jacobi
map is quasi-finite of degree 84 over an open subset
of the theta-divisor $\Theta\subset J(V_2)$.

In the index-1 case, several descriptions of the moduli
spaces $M_X(2;k,n)$ were obtained for the following threefolds:
the 3-dimensional
quartic \cite{IM-2}, the Fano threefold of degree
12 \cite{IM-3} and the one of degree 16 \cite{IR}. The vector bundles studied
in these three papers are related respectively to the half-canonical curves
of degree 15, elliptic quintics and elliptic sextics.
Kuznetsov in \cite{Ku-1},
\cite{Ku-2} used the moduli spaces 
associated to elliptic quintics on the 3-dimensional
cubic $V_3$ and the Fano threefold $X=X_{12}$ of degree 12
to construct semiorthogonal decompositions
of the derived categories of sheaves on these threefolds.

According to Mukai, any 
Fano threefold $X=X_{12}$ is a linear section of the
spinor tenfold in the projectivized half-spinor space
of Spin(10). The orthogonal linear section of the spinor
tenfold is a canonical \mbox{genus-7} curve $\Gamma$, and
the intermediate Jacobian $J(X)$ is isomorphic to the Jacobian
of $\Gamma$. It is proved in \cite{IM-3}
that $M_X(2;1,5)$ is isomorphic to~$\Gamma$. Kuznetsov
remarks that the last moduli space is fine and provides
a natural universal bundle on it.

Here we work on the same variety $X=X_{12}$, but consider
the moduli space $M_X(2;1,6)$. We prove that all the vector bundles
represented by points of $M_X(2;1,6)$ are obtained
by Serre's construction from reduced
sextics which deform to elliptic sextics (Proposition \ref{E-sextic}).
The main result
(Theorem \ref{singtheta} and Corollary \ref{MX-irred}) is the following:
for generic $X$, $M_X(2;1,6)$ is irreducible and
the Abel--Jacobi map sends it birationally onto
the singular locus $\Sing\Theta$ of the theta-divisor of $J(X)$.
Our construction provides no universal bundle on $M_X(2;1,6)$,
and it seems very likely that this moduli space is not fine.

Throughout the paper, we extensively
use the Iskovskikh--Prokhorov--Takeuchi
birational transformations
that can be obtained by a blowup with
center in a point $p$, a conic $q$ or a twisted rational cubic
$C_3^0$ followed by a flop
and a contraction of one divisor (Section~1).
The existence of such transformations is proved
in \cite{Tak}, \cite{IP} by techniques from Mori theory.
The principal idea is the following.
The anticanonical class $-K_{\tilde X}$
of the blowup ${\tilde X}$ of $X$
along one of the above centers is nef and big and
defines a small contraction of ${\tilde X}\lra W$ onto some
Fano 3-fold $W$ with terminal singularities.
By a result of Koll\'ar \cite{Kol-1},
there exists a flop ${\tilde X}\dasharrow {\tilde Y}$ over $W$.
The flop is a birational map,
biregular on the complement of finitely many flopping curves which are
exactly the curves
contracted to the singular points of $W$. The thus obtained variety
${\tilde Y}$ admits a birational contraction ${\tilde Y}\lra Y $
onto another Fano threefold $Y$ with Picard group $\ZZ$.
The composition $X\dasharrow \tilde{X}\dasharrow\tilde{Y}\lra Y$
is what we call an Iskovskikh--Prokhorov--Takeuchi transformation.

If one applies this construction to a conic $q$
in $X$, then the resulting birational map $\Psi_q$ (see Diagram \ref{fig-2})
ends up in the 3-dimensional quadric
$Q^3$, and the last blowdown in its decomposition
is the contraction of a divisor onto a curve $\Gamma_{10}^7\subset Q^3$
of genus 7 and degree 10. 
The curve $\Gamma_{10}^7$ is identified with the projection
of $\Gamma$, the orthogonal linear section associated to $X$,
from two points $u,v\in\Gamma$.
This allows us to parameterize the family of conics in $X$ by the
symmetric square $\Gamma^{(2)}$. Further, the rational normal quartics
$C_4^0$ in $X$ meeting $q$ at 2 points are transformed by $\Psi_q$ into
conics in $Q^3$ meeting $\Gamma_{10}^7$ in 4 points. If we denote the 4 points
$u_1,u_2,u_3,u_4$, then the divisor $u+v+\sum u_i$ on $\Gamma$
belongs to $W_6^1$. 
The Brill--Noether locus $W_6^1$ is nothing else
but the singular locus of the theta-divisor in $J(\Gamma )$,
and the Abel--Jacobi image of the degenerate elliptic sextic
$C_4^0+q$ is minus the class of $u+v+\sum u_i$. 
Any elliptic sextic in $X$ defines  a rank-2 vector
bundle $\EEE$ via Serre's construction ${\EuScript S}$. We show that the
fibers of ${\EuScript S}$ are the projective spaces
$\PP^3=\PP H^0(X,\EEE )$ and those of the Abel--Jacobi map on elliptic sextics
are finite unions of these $\PP^3$'s. Further, we verify that
the reducible sextics of type $C_4^0+q$ in
a generic fiber of the Abel--Jacobi map form an irreducible curve.
Hence the fiber of the Abel--Jacobi map is just one copy of $\PP^3$, which
implies the birationality part of the main result.

In order to handle degree-6 curves, we start with lines, conics,
then continue by rational normal quartics, each time constructing higher
degree curves as smoothings of the reducible one. Thus we prove auxiliary
results on the families of low degree curves which may be of interest
themselves. For example, we identify the curve $\tau (X)$ of lines in $X$
with tne Brill--Noether locus $W_5^1(\Gamma )$ and determine 
its genus $g_{\tau (X)}=43$ (Proposition \ref{lines}). We prove
that the surface of conics $\FFF (X)$
is isomorphic to $\Gamma^{(2)}$ (Proposition
\ref{conics}). This result was also obtained
by \cite{Ku-2} via a different approach using the Fourier--Mukai
transform $D^b(X)\lra D^b(\Gamma )$. 
It is curious to note that $\FFF (X)$ remains nonsingular
for {\em all} nonsingular $X$.

Proceeding to curves of higher degree, we
show that the families of rational normal cubics and quartics in $X$
are irreducible (Lemmas \ref{cubics}, \ref{irred-quartics}).
We prove that the family of degenerate
elliptic sextics of the form $C_4^0+q$ in $X$
is irreducible (Lemma \ref{good-sextics}). A standard monodromy
argument together with
the result of N.~Perrin \cite{P-2} on the irreducibility
of the family of elliptic curves of given degree on the spinor tenfold
$\Sigma$ allow us to deduce the irreducibility
of the family of elliptic sextics in~$X$ and that of the
moduli space $M_X(2;1,6)$.


On several occasions, we use the rigidity of the symmetric square of
$\Gamma$ in the following sense: $\Gamma^{(2)}$ has neither nontrivial
self-maps, nor maps
to a curve. Though the subject seems to be classical, we did not find
appropriate references and included the proof of the rigidity of
$\Gamma^{(2)}$ for a generic curve of genus $g\geq 5$ in the last section
(Proposition \ref{aut}).

\bigskip

{\sc Acknowledgements.} The authors thank Yu.~Prokhorov
and N.~Perrin for discussions.

\section{Preliminaries}

Let $\Sigma =\sten$ be the spinor tenfold in $\PP^{15}$.
It is a homogeneous space of the complex spin group $\Spin$,
the unique closed orbit of $\Spin$
in the projectivized half-spinor representation $\Spin :\PP^{15}\actson$.
It can be also interpreted as one of the two components
of the orthogonal Grassmannian
$G(4; Q)=\Sigma^+\sqcup\Sigma^-$ parametrizing the linear
subspaces $\PP^4$ of $\PP^9$ contained in a given smooth
8-dimensional quadric $Q=Q^8\subset\PP^9$.
See \cite{Mu-1}, \cite{RS} or Section 1 of \cite{IM-3} for more
details and for expicit equations of $\Sigma$.

The Fano threefold $\XT$ is a smooth 3-dimensional
linear section of $\Sigma$ by a subspace $\PP^8\subset\PP^{15}$.
We will also consider smooth linear sections of $\Sigma$
by linear subspaces $\PP^7$ and $\PP^6$, which are K3
surfaces, resp. canonical curves of degree 12. The Gauss
dual $\Sigma^\dual\subset\PP^{15\dual}$ of $\Sigma$
is naturally identified with $\Sigma$ via the so called
fundamental form on $\PP^{15}$, and to a linear section
$V=\PP^{7+k}\cap\Sigma$ for $k=-1,0$, resp. 1
we can associate the orthogonal linear section
$\check{V}=(\PP^{7+k})^\perp\cap\Sigma^\dual$. The orthogonal
linear section of a Fano 3-fold $\XT$ is a canonical genus-7
curve $\Gamma =\gseven$, and that of a K3 surface ($k=0$) is
another K3 surface. By \cite{Mu-1}, $\Gamma =\check{X}$ is not an arbitrary
smooth curve of genus 7, but a sufficiently generic one:
it has no $g_4^1$, neither $g_6^2$.

If we identify $\Sigma$ with $\Sigma^+\subset G(4; Q)$,
then $\Sigma^\dual$ is naturally identified with the
other component $\Sigma^-$ of $G(4; Q)$. Denote by $\PP^{15\pm}$
the half-spinor space spanned by $\Sigma^\pm$, so that
$\PP^{15+}=\PP^{15}$ and $\PP^{15-}=\PP^{15\dual}$.
For $c\in \Sigma^\pm$, introduce the following notation:

$\PP^4_c$, the linear subspace of $Q$ represented by $c$;

$\PP^{14}_c$, the tangent
hyperplane to $\Sigma^\mp$ in $\PP^{15\mp}$
represented by $c$;

$H_c$, the corresponding hyperplane section $\PP^{14}_c\cap\Sigma^\mp$;

$\ep (c)$, the sign of $c$, that is $\ep (c)\in\{ +, -\}$ and
$c\in\Sigma^{\ep (c)}$.

The following proposition lists some useful properties of
$\Sigma^\pm$.

\begin{proposition}\label{linsect}
The following assertions hold:

(i) For $c,d\in G(4; Q)$, $\ep (c)= \ep (d)$, that is $c, d$ lie in the same
component
of $G(4; Q)$, if and only if $\dim (\PP^4_c\cap\PP^4_d)\in \{ 0,2,4\}$.

(ii) For $c,d\in G(4; Q)$, $\ep (c)= -\ep (d)$, that is  $c,d$ belong to 
different components of $G(4; Q)$, if and only
if $\dim (\PP^4_c\cap \PP^4_d)
\in \{ -1,1,3\}$, where the negative dimension corresponds to the empty
set.

(iii) Let $c\in G(4; Q)$. Then $H_c =\{ a\in G(4; Q)\ 
\mid\ \dim (\PP^4_c\cap \PP^4_d)\in \{ 1,3\}\}
=\{ d\in \Sigma^{-\ep (c)}\ 
\mid\ \PP^4_c\cap \PP^4_d\not =\emptyset\}$.

(iv) The hyperplane $\PP^{14}_c$ is tangent to $\Sigma^{-\ep (c)}$ 
along a linear $4$-dimensional subspace $\PP^4\subset\Sigma^{-\ep (c)}$,
which we will denote by $\Pi^4_c$, and $\Pi^4_c=\{ d\in G(4; Q)\ 
\mid\ \dim (\PP^4_c\cap \PP^4_d)=3\}$. Any $3$-space
$\PP^3\subset Q$ determines in a unique way a pair
$\PP^4_c$, $\PP^4_d$ of $4$-subspaces of $Q$ containing
$\PP^3$, so $\Pi^4_c$ is naturally identified with the
dual of $\PP^4_d$.

(v)  $H_c$ is a cone whose
vertex ($=$ ridge) is $\Pi^4_c$ and whose base is the Grassmannian $G(2,5)$,
embedded in a standard way into $\PP^9\simeq (\Pi^4_c)^\perp$.
The linear projection with center $\Pi^4_c$ identifies the
open set $U_c=H_c\smallsetminus \Pi^4_c$ with 
the universal vector subbundle of $\CC^5\times G(2,5)$ of rank $3$.
\end{proposition}

\begin{proof} The assertions (i), (ii) are classical, see for example
\cite{Mu-1}. For a proof of (iii)--(v) see
\cite{IM-3}, Lemma 3.4.
\end{proof}

The families of lines and conics on the spinor tenfold are
easy to describe:

\begin{proposition}\label{lines-conics}
(i) Fix a plane $\PP^2$ contained in $Q=Q^8$. Then
$$
\ell_{\PP^2}^\pm=\{ c\in \Sigma^\pm\ \mid \ \PP^2\subset\PP^4_c\}
$$
is a line in $\Sigma^\pm$. Every line in $\Sigma^\pm$ is of this form.
The variety $\tau (\Sigma^\pm)$ is thus identified with the Grassmannian
$G(2;Q)$ parametrizing the planes $\PP^2$ contained in $Q$.

(ii) Fix a point $p\in Q$. Then
$$
Q^{6\pm}_p=\{ c\in \Sigma^\pm\ \mid \ p\in\PP^4_c\}
$$
is a nonsingular $6$-dimensional quadric contained in $\Sigma^\pm$.
Any conic $q$ in $\Sigma^\pm$ belongs to one of the following two
types: either $q$ lies in a plane $\PP^2$ contained in $Q$,
or there exist a unique point $p\in Q$ depending on $q$, and a plane
$\PP^2$ in the linear span $\PP^{7\pm}_p$ of $Q^{6\pm}_p$ such that
$q=Q^{6\pm}_p\cap\PP^2 $.

More generally, for any quadric $q^k$ of dimension $k=0,1,\ldots ,6$
contained in $\Sigma^\pm$, either its span $\PP^{k+1}$ is contained in
$\Sigma^\pm$, or there exists a unique point $ p\in Q$ such that
$\PP^{k+1}\subset \PP^{7\pm}_p$ and $q^k=\PP^{k+1}\cap Q^{6\pm}_p$.
\end{proposition}

\begin{proof} Assertion (i) is proved in \cite{RS}, Section 3.
For the part (ii), see \cite{Mu-1}, 1.14--1.15.
\end{proof}

We will often use the following property of the plane linear sections
of $\Sigma$, whose proof is obtained by a refinement of the proof
of Proposition 1.16 in \cite{Mu-1}:

\begin{lemma}\label{4-secant}
Let $\PP^2$ be a plane in $\PP^{15}$. If $\PP^2\cap\Sigma$
is finite, then $\length (\PP^2\cap\Sigma )\leq 3$.
\end{lemma}

Informally speaking, this means that $\Sigma$ has no $4$-secant $2$-planes.
As $\Sigma$ is an intersection of quadrics, any intersection
$\PP^2\cap\Sigma$ that contains a subscheme of
length 4 is either a line, or a line plus a point,
or a conic, or the whole plane $\PP^2$.

Let now $X=\XT$ be a smooth Fano threefold of degree 12.
We will describe the Iskovskikh--Prokhorov--Takeuchi
(\cite{IP}, \cite{Tak})
birational maps $\Phi_x$, $\Psi_q$, resp. $\Psi_{C_3^0}$ associated to
a point $x\in X$, a conic $q\subset X$, resp. a rational
normal cubic $C_3^0\subset X$ (Theorems 4.5.8, 4.4.11, 4.6.3 in \cite{IP};
see also Theorems 6.3 and 6.5 of \cite{IM-3} for
the first two). For the reader's convenience, we will briefly
remind their structure. Each of these maps is a composition
of three birational modifications: blowup of a point
or a curve in $X$, flop and blowdown of some divisor onto a curve.
The blowup gives a 3-fold $\tilde{X}$ with nef and big
anticanonical class and 2 exceptional divisors. The first one 
is that of the blowdown $\tilde{X}\lra X$. The contraction of the second one
provides a new 3-fold $Y$, but before the contraction, one has
to make a flop in finitely many irreducible curves
$C\subset\tilde{X}$ characterized by the condition $C\cdot K_{\tilde{X}}=0$.

Start by $\Phi_x$, the birational map 
associated to a generic point $x\in X$. It is a birational
isomorphism of $X$ onto $Y=Y_5$, 
the Del Pezzo variety of degree 5, that is a nonsingular
3-dimensional linear section $\PP^6\cap G(2,5)$ of the
Grassmannian in $\PP^9$. Its structure is described by Diagram 1:

\noindent
\begin{figure}[!h]
$$
\xymatrix{
E_X \ar @{^{(}-} @<-2.8pt> [r]
\ar @{-} @<2.4pt> [r]
\ar [d]
& \tilde{X}
\ar [d]_{\sigma_X}  
\ar@{-->}^{\text{flop}} [rr]
\ar [ddr]^(0.35){\tilde{\pi}}
&&\tilde{Y}
\ar @{-^{)}} @<-2.8pt> [r]
\ar @{-} @<2.2pt> [r]
\ar [d]^{\sigma_Y}
\ar [ddl]_(0.35){\tilde{\eta}}
& E_Y
\ar [d]   
\\
x \ar @{^{(}-} @<-2.6pt> [r]
\ar @{-} @<2.6pt> [r]
\ar@{-->} @/^1pc/ [u]
^{
\begin{array}{c}\mbox{\scriptsize blowup}\\ \text{\scriptsize of }
{\scriptstyle x}\end{array}
}
& X
\ar @{-}  [l]+<4pt,0pt>
\ar@{-->}^{\Phi_x} [rr]
\ar@{-->} [dr]_{\pi}
&& Y
\ar @{-^{)}} @<-2.6pt> [r]
\ar @{-} @<2.4pt> [r]
\ar@{-->} [dl]^{\eta}
& \Gamma
\ar@{-->} @/_1pc/ [u]
_{
\begin{array}{c}\mbox{\scriptsize blowup}\\ \text{\scriptsize of }
{\scriptstyle \Gamma}\end{array}
}\\
&&W&&
}
$$
\bigskip
\caption{The birational isomorphism \mbox{$\Phi_x:X\dasharrow Y_5$} $= $
$G(2,5) \cap {\bf P}^6$.}

\bigskip\label{fig-1}
\end{figure}

In the diagram, $\pi=\pi_{2x}$ is
the double projection from $x$, that is
the rational map $X\dasharrow \PP^4$ defined by the
linear system of hyperplanes in $\PP^8$ tangent to $X$ at $x$,
$\Gamma =\gseven$ is a canonical genus-7
curve contained in $Y$, and $\eta$ the projection
by the linear system of quadrics containing $\Gamma$. 
The map $\Phi_x$ is given by the incomplete linear
system $|{\mathcal O}_X(3-7x)|$ and the opposite
map $\Phi_x^{-1}$ by the linear system $|{\mathcal O}_Y(12-7\Gamma )|$.
The curve $\Gamma$ is isomorphic to the orthogonal linear
section $\Gamma=\check{X}$ of $\Sigma^-$, denoted by the same symbol.
Both projections $\pi,\eta$ are birational and end
up in the same singular quartic 3-fold $W\subset\PP^4$.
When lifted to $\tilde{X}$ and $\tilde{Y}$, they
become regular morphisms defined by the anticanonical
linear system: $\tilde{\pi}=\phi_{|-K_{\tilde{X}}|}$,
$\tilde{\eta}=\phi_{|-K_{\tilde{Y}}|}$. The essential
point in this diagram is that the flop $\tilde{X}\dasharrow\tilde{Y}$
is {\em a flop over} $W$, that is $\tilde{\pi}, \tilde{\eta}$
are small morphisms contracting the flopping curves
to isolated singular points of $W$, and these flopping curves
are the only indeterminacies of the flop. We showed in \cite{IM-3}
that for generic $X, x$, the flopping curves in $X$ are the 24
conics passing through $x$, and those in $Y$ are the 24
bisecant lines to $\Gamma$ contained in $Y$.

The map $\Psi_q :X\dasharrow Q^3$ of the second type is a 
birational isomorphism
from $X$ to a 3-dimensional quadric $Q^3\subset\PP^4$,  
associated to a generic conic $q\subset X$.  
It is given by
the linear system \mbox{$|\OOO_X(2-3q)|$}, and its inverse
$\Psi_q^{-1}$ by $|\OOO_{Q^3}(8-3\Gamma_q )|$. Its structure
is described by Diagram 2:  

\noindent
\begin{figure}[!h]
$$
\xymatrix{
E_X \ar @{^{(}-} @<-2.8pt> [r]
\ar @{-} @<2.4pt> [r]
\ar [d]
& \tilde{X}
\ar [d]_{\sigma_X}  
\ar@{-->}^{\text{flop}}_{\phi} [rr]
&&\tilde{Q}^3
\ar @{-^{)}} @<-2.8pt> [r]
\ar @{-} @<2.2pt> [r]
\ar [d]^{\sigma_{Q}}
& E_Q
\ar [d]   
\\
q \ar @{^{(}-} @<-2.6pt> [r]
\ar @{-} @<2.6pt> [r]
\ar@{-->} @/^1pc/ [u]
^{
\begin{array}{c}\mbox{\scriptsize blowup}\\ \text{\scriptsize of }
{\scriptstyle q}\end{array}
}
& X
\ar@{-->}^{\Psi_q} [rr]
&& Q^3
\ar @{-^{)}} @<-2.6pt> [r]
\ar @{-} @<2.4pt> [r]
& \Gamma_q
\ar@{-->} @/_1pc/ [u]
_{
\begin{array}{c}\mbox{\scriptsize blowup}\\ \text{\scriptsize of }
{\scriptstyle \Gamma_q}\end{array}
}
}
$$
\bigskip
\caption{The birational isomorphism $\Psi_q :X\dasharrow Q^3$.}

\bigskip\label{fig-2}
\end{figure}

In this diagram, $\Gamma_q\subset Q^3$ is a curve of degree 10 and genus 7,
isomorphic to the
orthogonal linear section $\Gamma$
associated to $X$ (see Corollary 5.12 in \cite{IM-3}).
It is not canonically embedded, for it has genus 7 and lies in $\PP^4$.
By the geometric Riemann--Roch Theorem, there is
a unique unordered pair of points $u, v\in\Gamma$ such that
$\OOO_{Q^3}(1)|_{\Gamma_q}\simeq \OOO_{\Gamma} (K-u-v)$,
where $K$ denotes the canonical class, and $\Gamma_q\subset\PP^4$
is the image of $\Gamma$ under projection
from the line $\bar{uv}$. We will denote it sometimes
$\Gamma_{u,v}$ in place of $\Gamma_q$. The flopping curves in
$X$ are the 14 lines meeting $q$, and those in $Q^3$ are the 14 trisecants
of $\Gamma_q$ contained in $Q^3$.

\noindent
\begin{figure}[!h]
$$
\xymatrix{
E_X \ar @{^{(}-} @<-2.8pt> [r]
\ar @{-} @<2.4pt> [r]
\ar [d]
& \tilde{X}
\ar [d]_{\sigma_X}  
\ar@{-->}^{\text{flop}} [rr]
&&\widetilde{{\PP^3}}
\ar @{-^{)}} @<-2.8pt> [r]
\ar @{-} @<2.2pt> [r]
\ar [d]^{\sigma_{\PP^3}}
& E_{\PP^3}
\ar [d]   
\\
{C_3^0} \ar @{^{(}-} @<-2.6pt> [r]
\ar @{-} @<2.6pt> [r]
\ar@{-->} @/^1pc/ [u]
^{
\begin{array}{c}\mbox{\scriptsize blowup}\\ \text{\scriptsize of }
{\scriptstyle C_3^0}\end{array}
}
& X
\ar@{-->}^{\textstyle{\Psi}_{C_3^0}} [rr]
&& {\PP^3}
\ar @{-^{)}} @<-2.6pt> [r]
\ar @{-} @<2.4pt> [r]
& {\Gamma_9^{7}}
\ar@{-->} @/_1pc/ [u]
_{
\begin{array}{c}\mbox{\scriptsize blowup}\\ \text{\scriptsize of }
{\scriptstyle \Gamma_9^{7}}\end{array}
}
}
$$
\bigskip
\caption{The birational isomorphism $\Psi_{C_3^0}:X\dasharrow \PP^3$.}

\bigskip\label{fig-3}
\end{figure}

The map $\Psi_{C_3^0}$ of the third type is a birational isomorphism
of $X$ onto $\PP^3$ and is described by Diagram 3.
In this diagram, $C_3^0$ is a sufficiently generic rational cubic
curve in $X$, and
$\Gamma_9^{7}\subset\PP^3$
is a nonsingular curve of degree $9$ and genus 7
which is a projection of the canonical curve $\Gamma =\check{X}$
from three points $u,v,w\in\Gamma$. The direct map 
$\Psi_{C_3^0}$ is given by the linear system $|\OOO_X(3-4C_3^0)|$
and its inverse by $|\OOO_{\PP^3}(15-4\Gamma_9^7 )|$.
The flopping curves in $X$ are the 21 lines meeting $C_3^0$,
and those in $\PP^3$ are the 21 quadrisecants to $\Gamma_9^7$.

\section{Lines and conics in $X_{12}$}

We will start the study of curves on $X$
with a description of the families of
lines and conics in terms of the orthogonal curve $\Gamma =\check{X}$.

\begin{proposition}\label{lines}
Let $X=X_{12}$ be any transversal linear section
$\PP^8\cap\Sigma$, $\Gamma =\check{X}$ its orthogonal
curve and $\tau (X)=\Hilb^{t+1}_X$ the Hilbert scheme
of lines in $X$, where a ``line'' is a subscheme of $X$
with Hilbert polynomial $P(t)=t+1$. Let $R(X)$ be the surface
swept by the lines in $X$: \ \ $R(X)=\bigcup\limits_{v\in
\tau (X)}\ell_v$.
Then the following statements hold.

(i) $\tau (X)$ is a connected locally complete intersection curve of arithmetic
genus $43$, isomorphic to the Brill--Noether locus $W_5^1(\Gamma )$.

(ii) If $X$ is generic, then $\tau (X)$
is nonsingular and every line $\ell\subset X$ has normal bundle
$\NNN_{\ell/X}\simeq \OOO_{\PP^1}(-1)\oplus \OOO_{\PP^1}$.

(iii) If $X$ is generic, then the generic line on
$X$ meets eight other lines and $R(X)\in |\OOO_X(7)|$.
\end{proposition}

\begin{proposition}\label{conics}
Under the hypotheses of the previous proposition, let
$\FFF (X)$ denote the Hilbert scheme $\Hilb^{2t+1}_X$
of conics on $X$ (the ``Fano surface" of $X$),
where a ``conic'' is a subscheme of $X$
with Hilbert polynomial $P(t)=2t+1$.
Then the following statements hold:

(i) A generic conic $q$ is
nonsingular and $\NNN_{q/X}\simeq \OOO_{\PP^1}\oplus \OOO_{\PP^1}$.

(ii) $\FFF (X)$ is isomorphic to $\Gamma^{(2)}$,
where $\Gamma^{(2)}$ denotes the symmetric square of $\Gamma$.

(iii) There are $24$ conics
passing through a generic point of $X$.
\end{proposition}

We will start by conics.

{\em Proof of Proposition \ref{conics}. }
For part (i), see \cite{IP}, Proposition 4.2.5, Remark 4.2.8
and Theorem 4.5.10. Part (iii) was proved in \cite{IM-3}, Theorem 6.3 (f).
We will now prove (ii).

We are going to construct an isomorphism
$\lambda:\FFF (X)\isoto\Gamma^{(2)}$. We will describe the
construction of $\lambda (q)$ for a closed point $q\in \FFF (X)$;
it is clear how one can extend it to $T$-points of $\FFF (X)$
for any scheme $T$.

Since $X$ is a linear section of $\Sigma=\Sigma^+$ and does not contain
planes, it does not contain conics of the first type
in the sense of Proposition \ref{lines-conics}. Hence to
any conic $q\subset X$ we can associate a unique point
$p=p(q)=\bigcap\limits_{x\in q}\PP^{4+}_x\in Q^8$, so that
$q=\PP^{2+}(q)\cap Q^{6+}_p$, where $\PP^{2+}(q)$ denotes
the linear span $\langle q\rangle$ of $q$. We can rewrite it as
$q=\PP^{8+}(X)\cap
\PP^{7+}_p\cap \Sigma^+\subset\PP^{15+} $, where
$\PP^{8+}(X)=\langle X\rangle$, $\PP^{7+}_p=\langle Q^{6+}_p\rangle$
and $\PP^{2+}(q)=\PP^{8+}(X)\cap\PP^{7+}_p $. If we now pass
to the orthogonal complements in $(\PP^{15+})^{\dual}=\PP^{15-}$,
we obtain: $$\PP^{12-}(q):=\PP^{2+}(q)^{\perp}=
\langle\PP^{8+}(X)^{\perp}, (\PP^{7+}_p)^{\perp}\rangle
= \langle\PP^{6-}(\Gamma), \PP^{7-}_p\rangle,$$
where $\PP^{6-}(\Gamma)=\langle\Gamma\rangle$. Thus
$\PP^{6-}(\Gamma), \PP^{7-}_p$ are not in general position,
but intersect in a line $\PP^1$. 
The triple intersection $\PP^{6-}(\Gamma)\cap \PP^{7-}_p\cap\Sigma^-$
can be seen as $(\PP^{6-}(\Gamma)\cap \PP^{7-}_p)\cap\Sigma^-=\PP^1\cap\Sigma^- $,
or $\PP^{6-}(\Gamma)\cap ( \PP^{7-}_p\cap\Sigma^-)=\PP^{6-}(\Gamma)\cap Q^{6-}_p$,
or else as 
$ \PP^{7-}_p\cap (\PP^{6-}(\Gamma)\cap\Sigma^-)=\PP^{7-}_p\cap \Gamma$.
Hence it is a subscheme of length 2 contained
in $\Gamma$, that is an element of $\Gamma^{(2)}$. We define:
$$\lambda (q):=\PP^{6-}(\Gamma)\cap \PP^{7-}_p\cap\Sigma^-\in\Gamma^{(2)}.$$

The inverse map is defined in exactly the same manner: by
Proposition \ref{lines-conics} (ii) for $k=0$, a subscheme
$\xi\subset\Gamma$ of length 2 is contained in a unique quadric
$Q^{6-}_p$ and we define:
$$\lambda^{-1}(\xi):=\PP^{8+}(X)\cap \PP^{7+}_p\cap \Sigma^+\in\FFF(X).$$
\hspace*\fill\square\smallskip

\begin{remark}
Alexander Kuznetsov \cite{Ku-2} proves the isomorphism
$\FFF (X)\simeq\Gamma^{(2)}$
in a more algebraic way: he shows that the Fourier--Mukai
transform associated to an appropriate universal rank-2 vector
bundle on $X\times M_X(2;1,5)$ sends the structure sheaf $\OOO_q$
of a conic $q\subset X$ to the sky-scraper sheaf $\OOO_\xi$ on $\Gamma
=M_X(2;1,5)$
for some $\xi\subset\Gamma$ of length~2.
\end{remark}

\begin{remark}
According to Mukai,
a generic K3 surface $S$ of degree 12 is a transversal linear
section of the spinor tenfold: $S=S+=\PP^7\cap \Sigma^+$.
Applying the same arguments as above with $\PP^7$ in place of $\PP^{8+}(X)$,
we obtain a non-isomorphic K3 surface $S^-=\PP^{7\perp}\cap\Sigma^-$
and an isomorphism 
$\lambda :\Hilb^2(S^+)\isoto\Hilb^2(S^-)$ (see also \cite{Mu-3}, Example 4).
\end{remark}

In our description of the birational transformation
$\Psi_q$ (see Diagram~\ref{fig-2}),
we associated a pair of points $u+v$ of $\Gamma$ to
a generic conic $q$. This gives
a rational map
$$
\FFF (X)\lra \Gamma^{(2)}\ ,\ \ q \mapsto u+v\ ,
$$
which we will temporarily denote by $f$. 

\begin{lemma}\label{lambda}
$\lambda =f$.
\end{lemma}

\begin{proof}
By Proposiition \ref{aut}, it suffices to prove
that $f$ is nonconstant. Then $f\circ \lambda^{-1}$
is a nonconstant rational self-map of $\Gamma^{(2)}$;
it is the identity  for generic $\Gamma$, and
by continuity, this is true for any nonsingular $\Gamma$.

Let
$u+v\in\Gamma^{(2)}$ be a generic degree-2 divisor.
Let $\Gamma_{u,v}$ be the curve of degree 10 in $\PP^4$ obtained as
the image of the canonical curve $\Gamma\subset\PP^{6-}(\Gamma)$
under the projection from the line $\bar{uv}$. It is contained
in a unique quadric $Q^3$. By \cite{Mu-2}, Theorem 8.1,
there is a Fano 3-fold $X'=X'_{12}$, defined as the non-abelian Brill--Noether
locus $M_\Gamma(2,K,3)$, and a smooth conic $q\subset X'$,
such that
$\Gamma_{u,v}$ together with its trisecants is the indeterminacy
locus of $\Psi_q^{-1}:Q^3\dasharrow X'$. Since a variety
$X_{12}$ is uniquely determined by its orthogonal curve $\Gamma$,
we have $X\simeq X'$, so $f$ is a dominant map, and this ends the proof.
\end{proof}


{\em Proof of Proposition \ref{lines}}.
By Shokurov's Theorem on the existence of lines, see 4.4.13 in \cite{IP},
and by ibid, Proposition 4.4.2, the scheme $\tau (X)$ is of
pure dimension 1, and the normal bundle of a line is of
type $(0,-1)$ if and only if this line is represented
by a nonsingular point of $\tau (X)$. So (ii)
is a consequence of (i) together with the smoothness of $W_5^1(\Gamma )$
for a generic curve $\Gamma$ of genus 7 (\cite{ACGH}, IV.4.4
and V.1.6).

Let us prove (i).
The easiest way to construct a map from $\tau (X)$ to $W_5^1(\Gamma )$
is by using either one of the birational maps $\Phi_x$ or $\Psi_q$
with generic $x$ or $q$. For example, let us do it for $\Psi_q$.

Let $\ell$ be a line in $X$ and $q$ a generic conic.
Then $\ell$ does not meet $q$ and $\tilde{\ell}=\Psi_q(\ell )$
is a conic. Recalculating the degree of $\ell$, equal to 1,
from the linear system that defines $\Psi_q^{-1}$, we see
that $\tilde{\ell}$ meets $\Gamma_q$ in a scheme $Z$ of length 5.
Denoting by angular brackets the linear span, we have
$\langle Z\rangle_{\PP^4} =\langle \tilde{\ell}\rangle_{\PP^4} =\PP^2$,
and if we pull back $Z$ to the canonical model
$\Gamma\subset\PP^6$ then we
will have $\langle Z+u+v\rangle_{\PP^6}=\PP^4$.
The latter linear span cannot be smaller than $\PP^4$, because $\Gamma$
has no $g_7^3$ (see \cite{Mu-1}). Hence $Z+u+v$ is an element of a $g_7^2$ and
$D_q=K-Z-u-v$
belongs to a $g_5^1$ on $\Gamma$. Thus we have constructed a
map

$$
\mu_q: \tau (X)\lra W_5^1 (\Gamma )\ ,\ \ \ell \mapsto
[D_q],
$$
where the brackets denote the class of a divisor in the Picard group.

Now let us verify that the inverse map $\mu_q^{-1}$ is well defined.
Take a point $z\in W_5^1 (\Gamma )$ representing a linear series
$g_5^1(z)$. Then $|K-z|$ is a $g_7^2$ and we have two cases:

Case A.
$|K-z -u-v|$ is a single effective divisor $Z$.

Case B. $|K-z -u-v|$ is a pencil $g_5^{1\prime}=\{Z(t)\}_{t\in\PP^1}$.

In the case A, projecting down to $\PP^4$, we get
a single conic $C_2^0(z)=\langle Z\rangle\cap Q^3$ meeting
$\Gamma_q$ in 5 points. 
Here we have two subcases: either $C_2^0(z)$ is irreducible,
or it is a reducible conic $\ell'_i\cup m$, where $\ell'_i$
is one of the trisecant lines of $\Gamma_q$, and $m$ is a
bisecant line.
When $C_2^0(z)$ is irreducible, we define
$\mu_q^{-1}$ at $z$ by $\mu_q^{-1}(z)=\Psi_q^{-1}(C_2^0(z))$,
where $\Psi_q^{-1}$ applied to a curve denotes the proper
transform of this curve.
In the other subcase, $\ell'_i$ is a flopping curve. It has no
proper transform in $X$, so $\mu_q^{-1}(z)$ should be determined
by considering a limit of the curves $\mu_q^{-1}(w)$ when $w\rar z$.
We use the following general observation concerning any flop $\phi$:
when a member $C_z$ of some algebraic family of curves $\{C_w\}_{w\in T}$
acquires an irreducible component which is a flopping curve, say $\ell'$,
then the limiting curve $D_z=\lim _{w\rar z}\phi^{-1}(C_w)$
of the flopped family $\{D_w\}_{w\in T}$
does not contain the flopping curve $\ell$, corresponding to
$\ell'$, and is the
proper transform of the remaining components of $C_z$:
$$
D_z =\lim _{w\rar z}\phi^{-1}(C_w)= \phi^{-1}(C_z\setminus \ell').
$$
Moreover, $D_z$ meets $\ell$ in this case.
Thus, when $C_2^0(z)=\ell'_i\cup m$, we should
put $\mu_q^{-1}(z)=\Psi_q^{-1}(m)$.

Now we will eliminate Case B. Assume that $|K-z -u-v|$ is a pencil.
Then we can associate to $z$ a pencil of conics
$C_2^0(z,t)=\langle Z(t)\rangle\cap Q^3$, and a pencil
of lines $\ell (z,t)$ in $X$, so that $\mu_q^{-1}$ is not defined
at $z$. The pair $u+v$ is
determined as the unique effective divisor in
$|K-z-Z(t)|$. 
On the other hand,
$u+v=\lambda (q)$. The generic point
of $W_{10}^4=K-\Gamma^{(2)}$ is not contained
in the image of the sum map
$W_5^1(\Gamma)\times W_5^1(\Gamma)\lra \Pic^{10}(\Gamma)$.
By dimension reasons, to see this, it is sufficient to
verify that for any $z\in W_5^1(\Gamma)$ there are finitely
many $w\in W_5^1(\Gamma)$ such that $|K-z-w|$ is effective.
This is stated in the following lemma.

\begin{lemma}\label{eight}
For the generic $z\in W_5^1(\Gamma)$ there are exactly $8$ distinct points
$w\in W_5^1(\Gamma)$ such that $|K-z-w|$ is effective.
\end{lemma}

\begin{proof}
The image $\bar{\Gamma}$ of $\Gamma$ under the map
given by the linear system $g_7^2=|K-z|$ is a plane septic
without triple points.
Hence $\bar{\Gamma}$ has exactly 8 double points, defining
8 linear subseries $g_5^1$ in the given $g_7^2$.
\end{proof}

Now we see that for a generic conic $q$, $\lambda (q)$
cannot be represented as the sum of two $g_5^1$'s,
hence Case B is impossible.

To compute the genus of $\tau (X)$,
we will use the  approach and the notation
from \S\ 8 of \cite{RS}.

Let $M$ be the base of the family of lines $\ell \subset \Sigma$
on the spinor 10-fold $\Sigma$.
By loc. cit.,
the incidence family
$$
G = \{ (x,L): x \in L \} \subset \Sigma \times M
$$
together with the natural projection
$\pr _1: G \rightarrow \Sigma$
is nothing else but the Grassmannization
$G  = G(3,\BBB) \rightarrow \Sigma$
of the universal subbundle
$\BBB \rightarrow \Sigma \subset G(5,10)$.

Let $h$ be the class of the hyperplane section of $\Sigma \subset {\PP}^{15}$
let $b_i = c_i(\BBB), i = 1,...,5$ be the Chern classes of $\BBB$,
and let $u_i = c_i(\UUU)$,  $i = 1,2,3$ be the Chern classes
of the universal subbundle $\UUU \subset \BBB_G$ on $G = G(3,\BBB)$;
in particular $-u_1 = -c_1(\UUU)$ is the class of the hyperplane section
of the Pl\"ucker embedding $M  \subset  G(3,10)$.
Then
$h^{10} = \deg\ \Sigma  = 12 \in {\QQ} = H^{20}(\Sigma , {\QQ})$,
\begin{equation}\label{a}
H^*(\Sigma , {\QQ}) \cong
{\QQ}[h,b_3]/(b_3^2+8b_3h^3+8h^6,6h^5b_3+7h^8),
\end{equation}
and the cohomology ring $H^*(G,{\QQ})$ is generated
as a $H^*(\Sigma , {\QQ})$-algebra by $u_1$ and $u_2$:

\begin{equation}\label{b}
H^*(G,{\QQ}) \cong
H^*(\Sigma, {\QQ})[u_1,u_2]/(f,g),
\end{equation}
where
$$
f  = h^4 - h^2u_2 - \frac{1}{2}u_2^2 - \frac{1}{2}b_3u_1
    + 2h^3u_1  -2hu_2u_1  + 3h^2u_1^2 -\frac{1}{2}u_2u_1^2
    + 2hu_1^3 +\frac{1}{2}u_1^4
$$
and
\begin{multline*}
g = b_3h^2 - \frac{1}{2}b_3u_2 + 2h^3u_2 - hu_2^2 + b_3hu_1
    + 3h^2u_2u_1 - 2u_2^2u_1 - \\ \frac{1}{2}b_3u_1^2 + 2h^2u_1^3
    + \frac{1}{2}u_2u_1^3  + 2hu_1^4 + \frac{1}{2}u_1^5.
\end{multline*}
In particular, the definition of the universal subbundle
$\UUU  \rightarrow G = G(3,\BBB)$ yields

\begin{multline}\label{c}
u_1^6h^{10} = (-u_1)^6h^{10} =  \deg G(3,5)\cdot \deg \Sigma
= \\
5\cdot 12 = 60 \in {\QQ} = H^{32}(G,{\QQ}).
\end{multline}

The second projection $\pr _2: G \rightarrow M$
is a projectivization of the rank-2 vector bundle
$\EEE = {\pr _2}_*\pr _1^*{\OOO}(h)$, and
$$
H^*(G,{\QQ}) \cong H^*(M,{\QQ})[h]/(h^2 - c_1h + c_2)
$$
where $c_1 , c_2$ 
are the Chern classes of $\EEE$. Thus $c_1 = -u_1$ and $c_2  = -h^2-u_1h$.
We have also $K_M = 6u_1$.

Since ${\tau}(X) \subset M$ is the common zero locus of $7$
general sections of $\EEE$, then
$[\tau  (X)] = c_2(\EEE)^7 =  (-h^2-u_1h)^7$ and
$K_{\tau (X)} = (K_M + 7c_1(\EEE))|_{\tau (X)} = -u_1|_{\tau (X)}$.
Therefore $\tau (X) \subset M \subset G(3,10)$ is a canonical
curve, and it remains to compute
the degree
$$
d = (-h^2-u_1h)^7(-u_1)h \in H^*(G,{\QQ})
$$
of $\tau (X)$ with respect to the Pl\"ucker
hyperplane class $-u_1$. This is done by reducing $d$
modulo the relations specified in (\ref{a}), (\ref{b}), (\ref{c}),
and the answer is $d=84$. 
Hence
$\tau (X) \subset G(3,10)$ is a canonical
curve of genus $g_{\tau (X)} = \frac{1}{2}d + 1 = 43$.

\medskip

To prove (iii), note that $\deg R(X)=\deg \tau (X)=84$,
hence $R (X)\sim 7H$.
For any line $\ell$, $\deg\NNN_{\ell/X}=-1$, so the contribution
of $\ell$ to the intersection number $\ell\cdot R(X)$ is $-1$,
hence $\ell$ meets $R(X)$ in eight isolated points counted with
multiplicities. As $X$ is generic, neither
of the lines on $X$ is a double curve of $R(X)$
and the multiplicity of a point of $R(X)$ equals the number
of lines passing through this point. Hence any line $\ell$
meets exactly 8 other lines.

\hspace*\fill\square\smallskip

\section{Abel--Jacobi map}\label{sect-AJ}

Let $X=X_{12}$ be any transversal linear section
$\PP^8\cap\Sigma$.
Let $J^d(X)$ denote the set of classes of algebraic 1-cycles
of degree $d$ in $X$ modulo rational equivalence.
It has a natural structure of
a principal homogeneous space under $J^0(X)$, and
according to \cite{BM}, $J^0(X)=J(X)$ is nothing else
but the intermediate Jacobian of $X$. Either of
the birational isomorphisms $\Phi_y,\Psi_q$ can be
used to identify
$J(X)$ with the Jacobian $J(\Gamma )=\Pic^0 (\Gamma )$. It is more convenient
to use $\Psi_q$. With the notation from the proof of Proposition \ref{lines},
the identification
goes as follows: $J(Q^3)=0$, and the passage from $Q^3$ to
$X$ consists in blowing up only one irrational curve $\Gamma_q$
followed by blowups of rational curves and their inverses.
By \cite{CG}, only the blowup with nonrational center
modifies the intermediate Jacobian, therefore
$J(X)\simeq J(Q^3)\times J(\Gamma_q )\simeq J(\Gamma )$.
This isomorphism is induced by the map $\Gamma_q\lra J^d (X)$,
$u\mapsto [\Psi_q^{-1}(u)]$, where
$d=\deg \Psi_q^{-1}(u)$. Here $\Psi_q^{-1}(u)$
is the image of the exceptional fiber $\sigma_Q^{-1}(u)\simeq\PP^1$
of $\sigma_Q$ over a point $u\in\Gamma_q$ under the map
$\sigma_X\circ\phi^{-1}$, where $\phi$ is the flop (see Diagram \ref{fig-2}).
It is irreducible
for generic $u$ and has a flopping
curve as one of its components for a finite set of values of $u$
corresponding to the points of intersection of trisecants with $\Gamma_q$.
According to Theorem 5.5 of \cite{IM-3}, the
curves $\Psi_q^{-1}(u)$ are
the rational cubics meeting $q$ twice. Applying the Abel--Jacobi functors
provides the desired isomorphism $a^1_q:\Pic^1(\Gamma )\isoto J^3(X)$.

As in loc. cit.,
we use the symbol ${\mathcal C}^g_d[k]_Z$
to denote the family of all the connected curves
of genus $g$ and degree $d$ meeting $k$
times a given subvariety $Z$ of a given variety $V$. More precisely, let
$Z\subset V$ be a nonsingular curve (resp. a point).
Then ${\mathcal C}^g_d[k]_Z$ is the closure in
the Chow variety of $V$ of the family of reduced connected curves
$C$ of degree $d$ such that 
length$\: (\OOO_X/(\III_C+\III_Z))=k$ (resp. $\mult_ZC=k$) and
$p_a(\tilde{C})=g$, where $\tilde{C}$ is the proper transform
of $C$ in the blowup of $Z$ in $V$.

We will summarize the above in the following lemma:

\begin{lemma}\label{aj}
Let $q$ be a generic conic in $X$. Then for any $k\in\ZZ$,
there is a natural isomorphism 
$$a^k_q:\Pic^k(\Gamma )\isoto J^{3k}(X),
\ \ \left[\sum n_iu_i\right]\mapsto 
\left[\sum n_i\Psi_q^{-1}(u_i)\right],$$
depending on $q$.

All the curves
$C\in{\mathcal C}^0_3[2]_q$, except for finitely many of them,
are irreducible and their images $\Psi_q(C)$ are
points of $\Gamma_q$. This yields a map $b_q :{\mathcal C}^0_3[2]_q
\lra \Pic^1(\Gamma )$. With the identification
$Pic^1(\Gamma )\isoto J^3(X)$ given by $a^1_q$, the map
$b_q$ is the Abel-Jacobi map
of the family ${\mathcal C}^0_3[2]_q$.
\end{lemma}

Now we will study the Abel--Jacoby map of more general families of curves
on $X$. We will use without mention the identification
of $J^k(X)$ and $\Pic^d(\Gamma )$, which is determined
by Lemma \ref{aj} uniquely modulo a constant translation.
Remark also that $J(X)=J(\tilde{Q^3})$ in a natural way.

\begin{lemma}\label{T}
Let $q$ be a generic conic in $X$. 
Let $T$ be the base of an irreducible family of curves on $X$ whose
generic member is a reduced curve which 
intersects neither $q$, nor any of the flopping curves of $\Psi_q$. 
Assume that
$\Psi_q$ transforms the family parameterized by
$T$ into a subfamily of $\CCC^g_d[k]_{\Gamma_q}$
on $Q^3$. For generic $C\in T$, denote by $Z_C$ or $Z^q_C$ the intersection
cycle $\Psi_q(C)\cap\Gamma_q$ considered as a degree-$k$
divisor on $\Gamma$. It can be defined by the formula
$Z_C=\sigma_{\Gamma *}(\tilde{C}\cdot E_Q )$, where
$\tilde{C}=\phi\sigma_X^{-1}(C)$ is the image of $C$ in $\tilde{Q}^3$
and $\sigma_\Gamma :E_Q\lra \Gamma_q$ is the restriction
of $\sigma_Q$.
Then the Abel--Jacobi map for the family
$T$ is given, up to a constant translation, by $C\mapsto -[Z_C]\in 
\Pic(\Gamma )$. 
\end{lemma}

\begin{proof}
As $J(Q^3)=0$, the Abel--Jacobi class
of the pullback of any family of curves on $Q^3$
is a point. Hence the class of $\sigma_Q^{-1}\Psi_q(C)$
in $J^\ccd (\tilde{Q}^3)$ is a constant, say $c$. If $Z_C=\sum n_iu_i$
($n_i\in\NN$, $u_i\in \Gamma$), then $[\sigma_Q^{-1}\Psi_q(C))]=
[C]+\sum n_i[\sigma_Q^{-1}(u_i)]$ and $[C]=c-\sum n_i[u_i]$,
as was to be proved.
\end{proof}

Now we will invoke the exceptional curves of $\sigma_X$.
By \cite{IM-3}, Theorem 5.5, their images in $Q^3$
are the elements of the family $\CCC_3^0[8]_{\Gamma_q}$.
Hence to each curve $\sigma_X^{-1}(x)$ with $x\in q$
we can associate a degree-8 divisor on $\Gamma$, defined
by $\sigma_Q\circ\phi(\sigma_X^{-1}(x))\cap\Gamma_q$. Its class
in $\Pic^8(\Gamma )$ does not depend on $x\in q$,
because $q$ is rational. Denote it by $d_8^q$.

\begin{lemma}\label{Tq}
In the hypotheses of the previous lemma, assume that
the generic curve $C_t$ of $T$ is of degree $d$ and does not meet $q$.
Let $C_0$ be a special member of $T$ such that the scheme-theoretic
intersection $C_0\cap q=M$ is of length $r$.  
Let $\tilde{C}_t$ be the pullback
of $C_t$ to $\tilde{X}$ for $t\neq 0$, and $\tilde{C}_0$ the limit of
$\tilde{C}_t$ as $t\rar 0$. Assume that
$C_0$ does not meet any of the flopping curves. Then the flop
$\phi$ is locally an isomorphism in the neighbourhood of
$\tilde{C}_0$
and all the nearby curves $\tilde{C}_t$, and the limit
of $[Z_C]$ when $t\rar 0$
is $\sigma_{\Gamma *}(\phi(\tilde{C}_0)\cdot E_Q )$. This
coincides with $[Z^q_{C_0}]+rd_8^q$ in the case when
neither of the components of $\tilde{C}_0$ is contracted
by $\sigma_{Q}$.
\end{lemma}

\begin{proof}
Let $M=\sum n_ix_i$. Then
$\tilde{C}_0=C_0'+\sum n_i\sigma_X^{-1}(x_i)$, where $C_0'$ is the proper
transform of $C_0$. The result follows by applying $\sigma_{Q*}$
to $\phi(\tilde{C}_t)\cdot E_Q$ as $t\rar 0$. 
\end{proof}

Remark that
$\sigma_{\Gamma *}(\sigma_q^{-1}(u)\cdot E_Q)=-u$, so the Abel--Jacobi
image of $\sigma_Q^{-1}(u)$ is $[u]$, which agrees with Lemma \ref{aj}.

In the proof of Propositions \ref{lines} and \ref{conics},
we introduced the maps
\mbox{$\mu_q :\tau (X)\lra W_5^1(\Gamma )$} and
$\lambda :\FFF (X)\lra \Gamma^{(2)}=W_2^0(\Gamma )$.
They can be considered as maps to $\Pic (\Gamma )$.

\begin{lemma}\label{mu}
The map $\mu=\mu_q$ does not depend on the choice of a generic conic $q$ and 
is, up to a constant translation, the Abel-Jacobi map of the family
of lines on $X$.
\end{lemma}

\begin{proof}
For generic $X$, $\mu_q$ is an isomorphism of two nonsingular
curves of genus $43$. A curve of genus $\geq 2$
has only finitely many automorphisms, 
hence $\mu_q$ does not depend on $q$ for generic $X$. 
As we saw in the proof of Proposition \ref{lines}, $\tau (X)$
remains a l.~c.~i. curve and is a
zero locus of a section of a vector bundle for all
nonsingular varieties $X$. Hence all of the components of $\tau (X)$
for the special (but still smooth) 3-folds $X$ are in the limit
of the family of curves $\tau (X)$ for nearby general 3-folds
$X$. Hence $\mu_q$ does not depend on $q$
by continuity on the special $X$, too.

The $\Psi_q$-image of a line $\ell$ not meeting $q$ is a conic
meeting $\Gamma_q$ in a degree-5 divisor $Z^q_\ell$, and
\begin{equation}\label{mu-q}
\mu_q(\ell )=K-\lambda (q)-[Z^q_\ell ].
\end{equation}
By Lemma \ref{T},
$\mu_q$ is, up to a constant translation, the Abel--Jacobi
map of the family of lines on $X$.
\end{proof}

In the following definition we generalize
the formula (\ref{mu-q}) to curves of any degree.

\begin{definition}\label{def-aj-can}
Let
$C\subset X$ be a curve of degree $d$, and $q$ a sufficiently generic
conic in $X$. This means that $q$ is not a component of $C$, $\Psi_q$
exists and $C$ does
not meet any of the flopping curves of $\Psi_q$. In this case the
scheme-theoretic inverse image $\tilde{C}=\sigma_X^*(C)$
is mapped isomorphically by the flop $\phi$ to a curve in $\tilde{Q}^3$. 
Let $\length (C\cap q)=r$ and $\Psi_q(C)\cap\Gamma=Z^q_C$.
Define 
\begin{multline}\label{for-aj-can}
AJ(C)=dK-d\lambda (q)-
\sigma_{\Gamma *}(\phi(\tilde{C})\cdot E_Q )=\\
d(K-\lambda (q))-rd_8^q-[Z^q_C]
\in \Pic^{5d}(\Gamma ).
\end{multline}
We call $AJ(C)$ the canonical Abel--Jacobi image
of $C$ in $\Pic^{5d}(\Gamma )$.
\end{definition}

Now we will determine the canonical Abel--Jacobi image of a conic.

\begin{lemma}\label{zqq}
For a generic pair of conics $q,q'$ on $X$,
$$[Z^q_{q'}]=K-2\lambda (q)+\lambda (q').$$
\end{lemma}

\begin{proof}
The $\Psi_q$-image of $q'$ in $Q^3$ is a rational quartic
$C_4^{q}(q') \subset Q^3$ intersecting ${\Gamma}_{q}$ in a divisor
$Z_{q'}^q$ of degree $10$.
>From the ideal sheaf sequence for $C_4^{q}(q') \subset Q^3$
we obtain
$$
h^0(\III_{C_4^{q}(q'), Q^3}(2))
\ge
h^0({\cal O}_{Q^3}(2)) - h^0({\cal O}_{C_4^{q}(q')}(2))
= 14-9=5.
$$

Therefore there exists a ${\PP}^4$-family of quadric sections
$S(t)$ of $Q^3$ through $C_4^{q}(q')$.
Each of these $S(t)$ intersects ${\Gamma}_{q}$ in a divisor
$D_{20}(t) \sim 2K - 2{\lambda}(q)$ of degree 20 such that
$D_{20}(t) = Z_{q'}^q + D_{10}(t)$ for an effective
divisor $D_{10}(t)$ of degree $10$ on $\Gamma$.
Therefore $h^0(D_{10}(t)) \ge 5$. Since $\deg \ D_{10}(t) = 10$
(and $\Gamma$ is non-hyperelliptic), we have $h^0(D_{10}(t)) \ge 5$
and $D_{10}(t) = K - D_2(t)$ for some divisor $D_2(t)$ of degree $2$.
Again, as $\Gamma$ is non-hyperelliptic, $D_2(t)$ does not depend
on $t \in {\PP}^4$.

Therefore $D_2(t) =D_2(q,q')$
depends only on $q$ and $q'$, and \
$Z_{q'}^q =  2H - D_{10}(t)$
=
$(2K -  2{\lambda}(q)) - (K - D_2(q,q'))$
=
$K - 2{\lambda}(q) + D_2(q,q')$.

If one regards $q$ as a fixed conic and $q'$ as a general one,
then the map $q'\mapsto -[Z_{q'}^q]$ is, up to translation,
the Abel--Jacobi map of the family of conics. It is obviously nonconstant.
Indeed, assume the contrary. Then any two conics are rationally
equivalent. Hence the sums $\ell +m$ of intersecting lines
are all rationally equivalent. This implies that $W_5^1(\Gamma )$
is hyperelliptic and the curve $F$ of pairs of intersecting lines
is a $g_2^1$ on it, hence $F$ is rational.
This is absurd, for $F\subset\tau (X)^{(2)}$ is mapped
injectively into $\FFF (X)$ and
$\FFF (X)\simeq\Gamma^{(2)}$ does not contain
rational curves. Therefore the Abel--Jacobi map of conics is nonconstant,
and hence the map $q'\mapsto D_2(q,q')$ is nonconstant as well.

Thus the composition of this map with $\lambda$ is a nonconstant
self-map of $\Gamma^{(2)}$. By Proposition \ref{aut},
it is the identity. Hence $D_2(q,q') = {\lambda}(q')$.
\end{proof}

\begin{corollary}\label{aj-conics}
The canonical Abel--Jacobi map $AJ|_{\FFF (X)}$ of
the family of conics on~$X$ is given by the formula
$$
AJ(q)=K-\lambda (q)\ \ \forall\ q\in\FFF (X).
$$
\end{corollary}

\begin{proposition}\label{aj-can}
The map $AJ$ defined by formula (\ref{for-aj-can})
does not depend on $q$, hence
$AJ$ induces a canonical isomorphism
$J^d(X)\isoto \Pic^{5d}(\Gamma )$ such that
$a^k_q\circ AJ$ is the translation  by a constant
depending only on $k,q,d$. For any two curves $C_1,C_2$
on $X$, we have
$$
AJ(C_1+C_2)=AJ(C_1)+AJ(C_2).
$$
\end{proposition}

\begin{proof}
By Lemma \ref{mu} and Corollary \ref{aj-conics}, the first statement
of the proposition is true for lines and conics. The statement
on the additivity
of $AJ$ is an immediate consequence of the definition, and we can
use it to extend the first statement from lines and conics to
curves of all degrees.

The Abel-Jacobi image of $\Gamma$ in $J(\Gamma)$ (defined
up to a translation) generates $J(\Gamma)$, hence the same is
true for the Abel-Jacobi image of $\Gamma^{(2)}$.
Hence the $AJ$-image of the family of conics generates $J(X)=J(\Gamma)$.
This means that
any algebraic 1-cycle on $X$ is rationally equivalent to
a linear combination of conics, and we are done.
\end{proof}

\begin{lemma}\label{d-eight}
For a generic conic $q \subset X$, the divisors
of the linear system $d_8^q$ on ${\Gamma}$, defined by
the intersections of the extremal rational cubics
$C_3 \in {\cal C}^0_3[8]_{\Gamma_q}$ with ${\Gamma}_q$,
belong to the linear system $|K - 2{\lambda}(q)|$.
\end{lemma}

\begin{proof}
We can assume $\Gamma$ (or $X$) generic;
the result for any $\Gamma$ will follow by continuity.
Consider the curve $D_q\subset \FFF (X)$ of conics
$q'$ in $X$ intersecting $q$, defined as the closure
of the set
$\{ q'\in\FFF (X)\mid \ q\cap q'\neq\emptyset ,\ 
\#(q\cap q')<\infty\}$. Let $q'\in D_q$. Then
$\Psi_q(q')$ is generically a bisecant line to $\Gamma_q$,
so that $Z_{q'}^q$ is a pair of points.
Using Corollary \ref{aj-conics},
Proposition \ref{aj-can} and Lemma \ref{Tq},
we can express the canonical Abel Jacobi image
of $q'$ in two different ways:
$$
AJ (q')=K-\lambda (q')=2K-2\lambda (q)-d_8^q-[Z_{q'}^q],
$$
where $Z_{q'}^q\in\Gamma^{(2)}$. Hence
$[Z_{q'}^q]=c+\lambda (q')$ for some constant $c=c(q)\in\Pic^0 (\Gamma )$
and for generic $q'\in D_q$.

Now extend this construction to the whole incidence 3-fold
$D$, the closure in $\FFF (X)\times\FFF (X)$
of the set $\{ (q,q')\mid q\cap q'\neq\emptyset ,\ 
\#(q\cap q')<\infty\}$.
Then we obtain the maps $h:D\lra \Gamma^{(2)}$, $(q,q')\mapsto Z_{q'}^q$, and
$c:\FFF (X)\lra J(\Gamma )$, $q\mapsto c(q)$, such that
$h(D)=\bigcup\limits_{q\in\FFF (X)}(c(q)+\lambda (D_q))\subset\Gamma^{(2)} $.
Assume that $c(q)\neq 0$ for some $q$. Then there is a one-parameter
family of distinct representations of $c(q)$ as the difference
$w(t)-z(t)$ of points $z(t),w(t)=z(t)+c(q)\in\Gamma^{(2)}$,
parameterized by $t\in D_q$.
Hence $w(t)+z(t')=w(t')+z(t)$ in $\Pic^4(\Gamma )$
for $t, t'$ moving in the same connected
component of $D_q$. This either implies the existence of
a linear series $g_4^1$ on $\Gamma$, or $D_q=u+\Gamma$, $c(q)=v-u$
for some $u,v\in\Gamma$. The first alternative is impossible, see \cite{Mu-1}.
The second one is also false. Indeed, the lines
spanned by the pairs $Z_{q'}^q$
for $q'\in D_q$ are secant lines
of $\Gamma_q$ contained in $Q$, but not all such secant lines
pass through a given point $v\in\Gamma_q$.
Hence $c(q)\equiv 0$ and we are done.
\end{proof}

\begin{corollary}
On the family $\CCC_3^0[2]_q$, the canonical Abel--Jacobi map
is given by
$$
AJ(C)=K+\lambda (q)+[\Psi_q (C)]\ \mbox{for generic}\ q\in\FFF (X)
\ \mbox{and}\ C\in\CCC_3^0[2]_q.
$$
\end{corollary}

\begin{proof}
In the notation of Proposition \ref{aj-can},
$\tilde{C}=\sigma_X^{-1}(x_1)+\sigma_X^{-1}(x_2)+\sigma_\Gamma^{-1}(u)$
for some $x_1,x_2\in q$, $u\in\Gamma_q$. Then 
$\sigma_{\Gamma *}(\tilde{C}\cdot E_Q )=2d_8^q-u$. The result now follows
from Proposition \ref{aj-can} and Lemma \ref{d-eight}.
\end{proof}

This still holds for a special cubic $C_3^0$ of the form
$q'_0+\ell$, where $q,q'_0,\ell$ intersect each other
with multiplicity 1. Then $\ell$ is a flopping line
of $\Psi_q$, and $q'_0$ is a special element of $D_q$
(notation from the proof of Lemma \ref{d-eight}).
The flopping
curve in $Q^3$  corresponding to $\ell$ is a trisecant
$\ell'$ to $\Gamma_q$, and if
$\ell'\cap \Gamma_q=u_1+u_2+u_3$, then the image of $q_0$
in $\tilde{Q}^3$ is the exceptional curve $\sigma_Q^{-1}(u_i)$
for one of the values of $i=1,2,3$, say $i=3$.
The limit of the proper transforms
of the curves $q'\in D_q$ as $q'\rar q'_0$ is the reducible curve
$\sigma_Q^{-1}(u_3)+\tilde{\ell}'$, so that $AJ(q_0')$ is given
by the same formula as above with $Z_{q'_0}^q=u_1+u_2$. This implies:

\begin{corollary}
If, in the above notation, $q,q'_0,\ell$ intersect each other
with multiplicity 1, then $AJ(\ell)=u_1+u_2+u_3+\lambda (q)$
and $AJ(q_0')=K-u_1-u_2$.
\end{corollary}

We can apply the results of this section to obtain some additional
information on lines, conics and the map $\Psi_q$. First, we
can characterize the curve of reducible conics in $\FFF (X)$.

\begin{lemma}\label{int-lines}
Let $\ell$, $m$ be two distinct lines in $X$.
Then $\ell\cap m\neq\empt$ if and only if $|K-\mu (\ell)-\mu (m)|$
is nonempty. In this case $\ell\cup m$ is a reducible conic and
$\lambda(\ell\cup m)= K-\mu (\ell)-\mu (m)$.
\end{lemma}

\begin{proof}
This follows immediately from the existence of the
canonical Abel--Jacobi map $AJ$ such that
$AJ(\ell\cup m)=AJ(\ell)+AJ(m)$ and from Lemma \ref{mu}
and Corollary \ref{aj-conics}.
\end{proof}

The next lemma answers the question, which lines should be
considered as lines ``intersecting themselves''.

\begin{lemma}\label{double-lines}
Let $\ell$ be a line in $X$. Then there is a double
structure on $\ell$ making it a conic in $X$ if and
only if $\ell$ is a singular point of $\tau (X)$.
\end{lemma}

\begin{proof}
Assume that the normal sheaf of $\ell$ is
$\OOO_{\PP^1}\oplus \OOO_{\PP^1}(-1)$. Let $C$ be a
{\em plane} double
structure on a line $\ell$ , that is,
a double structure embeddable into $\PP^2$.
Any Gorenstein doubling $\ell$ is given
by Ferrand's construction \cite{F}, \cite{BanF} and is associated to
a surjective morphism of $\OOO_\ell$-modules
$\NNN_{\ell /X}^\dual\lra\LLL$, where $\LLL$ is
some invertible sheaf on $\ell$. The kernel of the surjection
can be represented in the form $\JJJ /\III_{\ell}^2$
for an ideal sheaf $\JJJ \subset\OOO_X$, and
this ideal sheaf defines the Ferrand's double structure $C$
on $\ell$: $\JJJ =\III_C$. 
The dualizing sheaf of Ferrand's double structure satisfies
$\omega_C|_{\ell}\simeq \omega_{\ell}\otimes\LLL^{-1}$.
Applying this to our situation, we see that $\LLL\simeq\OOO_\ell (k)$
for some $k\geq 0$, hence $\omega_C|_{\ell}\simeq \OOO_\ell  (-2-k)$,
which contradicts the property $\omega_C|_{\ell}\simeq \OOO_\ell  (-1)$
verified for a {\em plane} doubling of $\ell$.

For a line $\ell$ with normal sheaf $\OOO_{\PP^1}(1)\oplus \OOO_{\PP^1}(-2)$,
the surjection $\NNN_{\ell /X}^\dual\lra \OOO_{\PP^1}(2)$ defines a unique
plane double structure on $\ell$.
\end{proof}

\begin{corollary}\label{8-7}
Let $\ell$ be a generic line on $X$. Then there are exactly
$8$ distinct lines $\ell_i$ such that $\ell+\ell_i$ is a conic.
They satisfy the condition $K-\mu (\ell)-\mu (\ell_i)\in \Gamma^{(2)}$.

If the normal bundle of $\ell$ is of type $(0,-1)$, then
the lines $\ell_i$ meet $\ell$ and are different from $\ell$.

If the normal bundle of $\ell$ is of type $(1,-2)$, then
$K-2\mu (\ell)\in \Gamma^{(2)}$. In this case, only
one of the $\ell_i$ coincides
with $\ell$ and the $7$ others are distinct and different from $\ell$.
\end{corollary}

\begin{proof}
This follows from Lemmas \ref{eight} and \ref{int-lines}.
\end{proof}

%
%
%
%
%

\section{Rational normal curves in $X$}

Let $X=X_{12}=\PP^8\cap\Sigma$ be a Fano 3-dimensional linear section 
of the spinor tenfold $\Sigma$ and $\Gamma =\check{X}$ its orthogonal
curve.
We will use the symbol $\CCC_d^g(X)$,
or simply $\CCC_d^g$,
to denote some families of degree-$d$ curves of genus $g$ in $X$,
whose precise definitions will be given in the context, and
$C_d^g$ to denote a member of such a family.
A {\em rational normal curve} of degree $d$ in $X$
is an irreducible nonsingular curve $C$ in $X$ such that $\deg C=d$
and $\dim\langle C\rangle=d$. Let $\CCC_d^0$
be the family of rational normal curves of degree $d$ in $X$.

\begin{lemma}\label{cubics}
The family $\CCC_3^0(X)$
of rational normal cubics in $X$ is irreducible, $3$-dimensional
and is birational to the symmetric cube $\Gamma^{(3)}$ of the curve $\Gamma$.
The normal bundle of a generic $C_3^0\in \CCC_3^0(X)$
is $\OOO_{\PP^1}\oplus\OOO_{\PP^1}(1)$.
\end{lemma}

\begin{proof}
The family of rational normal cubics in $X$ was studied in
\cite{IP}, 4.6.1--4.6.4. The authors determined the normal bundle and
constructed the birational
transformations $\Psi_{C_3^0}$ associated
to sufficiently general rational normal cubics $C_3^0\in \CCC_3^0(X)$.
We will explain how the existence of these transformations
implies the irreducibility of $\CCC_3^0(X)$.
Theorem 4.6.4 in loc. cit. provides the inverse
construction, which permits
to reconstruct $C_3^0\subset X$ starting from any sufficiently general
$\Gamma_9^{7}\subset\PP^3$.
Recall that $\Gamma_9^{7}$ is a projection of $\Gamma\subset\PP^6$
from a unique triple of points $u,v,w\in\Gamma$.
Hence the open part of $\CCC_3^0(X)$ consisting of those cubics
$C_3^0$ for which the map
$\Psi_{C_3^0}$ exists is birational to the symmetric cube $\Gamma^{(3)}$
and is irreducible.
\end{proof}

\begin{lemma}\label{quartics}
Let $C$ be a connected Cohen--Macaulay curve of degree $4$ in $X$.
Then the following assertions hold:

(i) $\dim\langle C\rangle=4$.

(ii) If $C$ is reduced and irreducible, then it is a rational normal quartic.

(iii) If $C$ is the union of two conics $q\cup q'$ such that
$q\cap q'\neq\emptyset$ and $\# (q\cap q')<\infty$, then
$q$, $q'$ meet each other quasitransversely at a single point.

(iv) $C$ has no singular points of multiplicity $\geq 3$ and $p_a(C)=0$.
\end{lemma}

\begin{proof}
All the assertions are easy consequences of the fact that
$\langle C\rangle =\PP^4$.
The latter follows from
the non-existence of 2-planes that are 4-secant to $X$. 
Indeed, assume that $\langle C\rangle =\PP^3$.
Then for any plane $\PP^2$ in this $\PP^3$, $\PP^2\cap X$
contains at least the 4 points of $\PP^2\cap C$ (counted with
multiplicities). 
By Lemma \ref{4-secant},
we obtain a three-dimensional family
of conics or lines in $X$, which is absurd. 
If $\langle C\rangle =\PP^2$, then $\langle C\rangle\cap X $
is an intersection of quadrics, hence coincides with $\langle C\rangle =\PP^2$.
This is absurd, as $X$ does not contain planes.
Hence $\langle C\rangle =\PP^4$.
\end{proof}

\begin{lemma}\label{irred-quartics}
Let $X$ be generic. Then the family $\CCC_4^0(X)$
of rational normal quartics in $X$ is irreducible.
\end{lemma}

\begin{proof}
Consider the family $I$ of all pairs $(C_4^0,X)$,
where $X$ is a Fano 3-fold section of the spinor 10-fold $\Sigma$
and $C_4^0\in\CCC_4^0(X)$. It has two natural projections
$p:I\lra \CCC_4^0(\Sigma )$ and $q:I\lra G(9,16)$, \mbox{$q:X 
\mapsto \langle X\rangle=$}$\PP^8\subset\PP^{15}$,
where $\CCC_4^0(\Sigma )$
is the family of rational normal quartics in $\Sigma$.
A nonempty fiber $q^{-1}(u)$ is the family $\CCC_4^0(X_u)$,
where $X_u=\PP^8_u\cap\Sigma$, and $p^{-1}(C_4^0)$
is an open subset of the Grassmannian $G(4,11)$
parametrizing the subspaces $\PP^8\subset\PP^{15}$
which contain  $\PP^4=\langle C_4^0\rangle$. By a standard
monodromy argument, the irreducibility of the generic
fiber $q^{-1}(u)$ will follow from the following two facts:
(1) $I$ is irreducible; (2) simultaneously for all sufficiently general $u$,
 one can choose in the fiber $q^{-1}(u)$
one distinguished irreducible component depending rationally on $u$.
As the fibers of $p$ are irreducible, the first fact
is equivalent to the irreducibility of $\CCC_4^0(\Sigma )$.
The latter follows from \cite{P-1}, where the author
proves that the Hilbert scheme $\Hilb_\Sigma^\alpha$ of
irreducible nonsingular rational curves of class $\alpha$
in a complex projective homogeneous manifold $\Sigma$
is smooth and irreducible when $\dim\Sigma \geq 3$ and
$\alpha$ is strictly positive.
The last condition holds in our situation, because
$\Pic \Sigma\simeq\ZZ$.

Now we will produce a distinguished
component $\CCC_4^{0*}$ of $\CCC_4^0(X)$ for a fixed $X$.
Let $C_3^0$ be a generic rational normal cubic in $X$. It intersects
the surface $R(X)$ swept by the lines in $X$ at a finite number of points.
Hence there is at least one line $\ell$ in $X$ meeting $C_3^0$.
Such a line cannot intersect $C_3^0$ in a scheme of length $\geq 2$,
for then the quartic $C_3^0\cup\ell$ will span $\PP^3$,
which is impossible by Lemma \ref{quartics}.
Therefore the family $\CCC_{3,1}^0$ of reducible
quartics $C_4^0=C_3^0\cup\ell$, where
$C_3^0\in\CCC_3^0$, $\ell$ is a line and $\length (C_3^0\cap\ell )=1$,
is a finite cover of $\CCC_3^0$. It is 3-dimensional. By
the standard normal bundle sequence for a reducible nodal curve,
$\chi (\NNN_{C_4^0/X})=4$, so $\dim_{[C_4^0]}\Hilb_X\geq 4$ and
hence $C_4^0$ can be deformed into a smooth rational normal quartic.
We define $\CCC_4^{0*}$ to be the component containing the smoothings
of curves from $\CCC_{3,1}^0$, but for this we need to prove the
irreducibility of $\CCC_{3,1}^0$.

Let $C_3^0\in \CCC_3^0$ be sufficiently generic. Then
the lines $\ell$ such that $C_3^0\cup\ell\in \CCC_{3,1}^0$
are the flopping curves of $\Psi_{C_3^0}$ (see Diagram 3 of 
Section 1). Hence they are
in a bijective correspondence with the quadrisecants of $\Gamma_9^7$.
Let $u+v+w\in\Gamma^{(3)}$ be the triple of points of $\Gamma$
associated to $\Gamma_9^7$. Let $L$ be a quadrisecant of $\Gamma_9^7$
and $L\cap\Gamma_9^7=u_1+u_2+u_3+u_4$. Then the span of the
divisor $D=u_1+u_2+u_3+u_4+u+v+w$ in $\PP^6$ is $\PP^4$, hence
$D$ belongs to a linear series $g_7^2$. Let us denote by
$G_d^r$ the subset of $\Gamma^{(r)}$ which is the union of all the
linear series $g_7^2$. 
As a generic $\Gamma$ has no $g_7^3$, the natural map $\pi :G_7^2\lra W_7^2$
is a $\PP^2$-bundle over the smooth curve
$W_7^2\simeq W_5^1$ and the quadrisecants of $\Gamma_9^7$
are in a bijective correspondence with the elements of the subset
$\{ D\in G_7^2\mid D-u-v-w\ \mbox{is effective}\}$.

Let $I^{(k)}=\{ (F,D)\in \Gamma^{(k)}\times G_7^2\mid 
D-F\ \mbox{is effective}\}$ ($1\leq k\leq 7$),
and let
$q_k:I^{(k)}\lra G_7^2$ be
the natural projection. We have identified a dense open subset of
$\CCC_{3,1}^0$ with that of $I^{(3)}$. So we have to show that
$I^{(3)}$ is irreducible. The map $q_3:I^{(3)}\lra G_7^2$
is a 35-sheeted covering obtained by applying the relative
symmetric cube to the 7-sheeted covering $q_1$.
Hence it suffices to prove that the
monodromy group $M$ permuting the sheets of $q_1$ is the whole of $S_7$.
This follows from two facts: (a)~$M$~is transitive, that is
$I^{(1)}$ is irreducible, and (b)~$M$~is generated by transpositions.

To verify (a), restrict $q_1$ to the fiber $\PP^2$ of $\pi$
over a general $g_7^2\in W_7^2$. An orbit of length $k$ of $M$ gives
rise to a $k$-valued multisection of $q_1|_{q_1^{-1}(\PP^2)}$,
or equivalently, to a map $\PP^2\lra \Gamma^{(k)}$. But $\Gamma^{(k)}$
does not contain rational curves for $k<5$, since $\Gamma$
has no linear series of degree $k<5$. If we assume
that $M$ is not transitive, then there is an orbit of length
$k<4$ and the above map $\PP^2\lra \Gamma^{(k)}$ is constant,
which immediately leads to a contradiction. Hence $M$ is transitive.

To verify (b), one can show that the ramification of
$q_1|_{q_1^{-1}(\PP^2)}$ is simple in codimension 1.
This follows from the observation that all the divisors
from the linear series $g_7^2$ are obtained as the intersections $L\cap\Gamma_0$,
where $\Gamma_0\subset\PP^2$ is the image of $\Gamma$ under the map
given by the $g_7^2$ and $L$ runs over the lines in $\PP^2$.
The ramification points of $q_1$ correspond to the points
of tangency of $L$ to $\Gamma_0$, and the ramification is simple
when $L$ is a simple tangent to $\Gamma_0$. But for $g_7^2$
generic, $\Gamma_0$ is a nodal septic of genus 7 having only finitely many
flexes or bitangents. Hence the ramification of $q_1$ is simple
in codimension 1.
\end{proof}

\begin{lemma}\label{quartics-2}
Let $X$ be generic. Then the family $\CCC_4^0(X)$
of rational normal quartics in $X$ is $4$-dimensional
and the normal
bundle of a generic quartic $C_4^0\in\CCC_4^{0}(X)$
is either $2\OOO_{\PP^1}(1)$ 
or $\OOO_{\PP^1}\oplus\OOO_{\PP^1}(2)$.
\end{lemma}

\begin{proof}
Take a generic pair of intersecting conics $q\cup q'$.
Both $q$ and $q'$ have normal bundle $2\OOO$.
The strong smoothability of $q\cup q'$ is proved
by a standard application of the Hartshorne--Hirschowitz techniques
\cite{HH}, so $q\cup q'$ is represented by a smooth point
of the closure $\bar{\CCC}_4^{0}$ 
of $\CCC_4^{0}$ in $\Hilb_X$ and 
$\dim \CCC_4^{0}=\chi (\NNN_{q\cup q'/X})=4$. For an example
of such argument see Lemma 1.2 of \cite{MT-2}.

>From the semicontinuity
of $h^1(\NNN_{C/X})$ for $C\in\bar{\CCC}_4^{0}$ and the fact that
$h^1(\NNN_{q\cup q'/X}(-x))=0$ for a point $x\in q\smallsetminus q'$,
we deduce that $h^1(\NNN_{C/X}(-x))=0$ for generic $C\in \CCC_4^{0}$
and $x\in C$.
This implies the assertion on the normal bundle.
\end{proof}

\begin{lemma}\label{quartics-normal}
Let $X$ and $C_4^0\subset X$ be generic. Then 
$\NNN_{C_4^0/X}\simeq 2\OOO_{\PP^1}(1)$.
\end{lemma}

\begin{proof}
Assume that $\NNN_{C/X}\simeq\OOO_{\PP^1}\oplus\OOO_{\PP^1}(2)$
for generic $C=C_4^0$.
Let $p\in C$ be a generic point.
Let $H(p)\subset \Hilb_X$ be the the closure 
of the family of rational normal quartics in $X$
passing through $p$.
It can be identified with a closed subscheme of $\Hilb_{\tilde{X}}$,
where $\tilde{X}$ is the blowup of $p$ in $X$. Let $\tilde{C}$
be the proper transform of $C$ in $\tilde{X}$. We have
$\NNN_{\tilde{C}/\tilde{X}}\simeq\OOO_{\PP^1}(-1)\oplus\OOO_{\PP^1}(1)$
and the tangent space to $H(p)$ at $[C]$ is identified with
$H^0(\NNN_{\tilde{C}/\tilde{X}})$. Since
$H^1(\NNN_{\tilde{C}/\tilde{X}})=0$, $H(p)$ is smooth at $[C]$
of dimension 2, and there is a unique component
$H(C,p)$ of $H(p)$ containing $C$. By our assumption, the proper transform
$\tilde{C'}$ of a generic quartic $C'$
in $H(C,p)$ has the same normal sheaf. Let $F$ be the
universal family of curves $\tilde{C'}$ over $H(C,p)$
and $\pi :F\lra \tilde{X}$ the natural map. For generic
$C'$, considered as a fiber of $F$ over the point $[C']\in H(C,p)$,
we have $\NNN_{C'/F}\simeq 2\OOO$, and for its image
$\tilde{C'}=\pi (C')$ in $\tilde{X}$, \
$\NNN_{\tilde{C'}/\tilde{X}}\simeq\OOO(-1)\oplus\OOO(1)$.
As any map from $2\OOO$ to $\OOO(-1)\oplus\OOO(1)$ has
its image in the second factor $\OOO(1)$,
the differential of $\pi$ is degenerate
at the points of $C'$, hence also at the generic point of $F$.
By the Sard theorem, $\pi (F)$ is a surface. Hence all the rational normal
quartics in $H(C,p)$ sweep out a surface in $X$, say $S(C,p)$.

Let $p'\neq p$ be another generic point of $C$. Then there is a 1-dimensional
family of rational normal quartics in $S(p)$ passing through both $p,p'$.
The curves of this 1-dimensional family cover an open set
of $S(C,p)$ and of $S(C,p')$. This implies that $S(C,p')=S(C,p)$.
We can also replace $C$ by a generic
curve $C'$ in $H(C,p)$, then take generic $p''\neq p'$ in $C'$ and see
that $S(C,p)=S(C',p')=S(C',p'')$. This implies, in particular, that $S(C,p)$
is generically smooth at $p$ and
that $S(C,p)$ contains a 3-dimensional family of
rational normal quartics. The 4-dimensional 
family of rational normal
quartics in $X$ is thus rationally fibered
over some irreducible curve $B$ into
3-dimensional families $H_t$, $t\in B$, such that the curves
parameterized by $H_t$ for generic $t$ cover a surface $S_t\subset X$
($S_t=S(C,p)$ for some $p\in C$, $C\in H_t$). The existence of
a three-dimensional covering family of rational curves implies
the rationality of $S_t$ for generic $t$.

Let $S=S_t$ for generic $t\in B$ and $T$ the minimal desingularization
of $S$. We have already seen that there are 1-dimensional families
of rational normal quartics passing through two generic points
$p,p'$ of $S$. By the ``bend and break argument"
(\cite{Kol-2}, Corollary II.5.6), there is a reducible member
in every such pencil. Thus $S$ is covered by conics.
As it is a rational surface, there is a linear pencil of
conics on $S$. This contradicts the non-existence of
rational curves in the symmetric square of $\Gamma$
and Proposition \ref{conics}.
\end{proof}

\begin{lemma}\label{quartics-span}
Let $X$ and $C_4^0\subset X$ be generic. Then 
$C_4^0$ is the scheme-theoretic intersection
$\PP^4\cap X$, where $\PP^4=\langle C_4^0\rangle$ is the linear
span of~$C_4^0$.
\end{lemma}

\begin{proof}
If $\PP^4\cap X$ contains a point $p$ of the secant 3-fold of $C_4^0$,
then $\PP^4\cap X$ contains also the secant line $\ell$ of $C_4^0$
passing through $p$, because $X$ is an intersection of quadrics.
But a generic $C=C_4^0$ in $X$ has no secant lines contained in $X$.
Indeed, if $p,p'$ are the points of $\ell\cap C$ for a line $\ell\subset X$,
then the infinitesimal deformations of $C$ with fixed point $p$
are given by $H^0(\NNN_{\tilde{C}/\tilde{X}})$ in the notation of
Lemma \ref{quartics-normal}. We have 
$\NNN_{\tilde{C}/\tilde{X}}\simeq 2\OOO$, so the infinitesimal
deformations lift to algebraic ones and $H^0(\NNN_{\tilde{C}/\tilde{X}})$
generates $\NNN_{\tilde{C}/\tilde{X}}$ at $p'$, hence $C$
can be moved off $\ell$ near $p'$ inside the family of quartics
passing through $p$. The line $\ell$ cannot deform with
$C$, for $C$ meets the surface swept by lines in a finite number
of points, and there are only finitely many lines through $p$
(Proposition 4.2.2, (iv) of \cite{IP}).

Assume now that $\PP^4\cap X$ contains a point $p$ not on
the secant threefold of $C$. Then there is a 1-dimensional family
of 3-secant planes $\PP^2_t$ to $C$ through $p$, parameterized
by the points $t$ of some curve $B$. These planes are
4-secant to $X$, hence, by Lemma \ref{4-secant},
they contain conics  $q_t$ lying in $X$ and passing through
$p$. All these conics lie in $\PP^4=\langle C\rangle$
and sweep out a surface in $X$. But $\Pic X\simeq \ZZ$,
so the linear span of a surface
in $X$ is at least $\PP^6$. The obtained contradiction
proves that $\PP^4\cap X=C$ set-theoretically.

%

Assume now that $\PP^4\cap X$
has an embedded component supported at $z\in C$. Then there is a line
$L\neq T_pC$ in $\PP^4$ passing through $z$ and tangent to $X$. 
Choose any $p\in L\smallsetminus\{ z\}$. Then there exists a
3-secant $\PP^2$ to $C$ passing through $p$ and $z$. It is
4-secant to $X$, as the intersection of $\PP^2$ with $X$ at $z$ is
multiple. Hence $\PP^2\cap X$ is a conic $q$ passing through
$z$ in the direction of $L$. Then $ \PP^4\cap X$ contains
a point $p'$ of $q$ which does not lie in the secant variety
of $C$, which contradicts to what we have proved.
\end{proof}

\section{Elliptic sextics in $X$}

An {\em elliptic sextic} in $X$ is a nonsingular
irreducible curve $C\subset X$ of genus $1$ and of
degree 6. We will also deal with degenerate ``elliptic"
sextics, which we will call just {\em quasi-elliptic sextics}.
A quasi-elliptic sextic is a locally
complete intersection curve $C$ of degree 6 in $X$,
such that $h^0(\OOO_C)=1$ and the canonical sheaf of $C$ is trivial:
$\omega_C=\OOO_C$. A reduced quasi-elliptic sextic will be
called a {\em good sextic}.

\begin{lemma}\label{good-sextics}
Let $q$ be a generic conic on $X$. Then $X$ contains a
$2$-dimensional family of good sextics of the form $C_4^0\cup q$,
such that $C_4^0$ is a rational normal quartic
and $\length (C_4^0\cap q)=2$, that is
$C_4^0, q$ meet each other quasitransversely in $2$ distinct points or
are mutually tangent at a single point. For a generic sextic
of this form, $C_4^0\cap q$ is a pair of distinct points.

If we let $q$ vary, then the family $\CCC_{4,2}^1$
of good sextics of type $C_4^0\cup q$ is
irreducible and $4$-dimensional. 
 
For any good sextic $C$ in $X$, $\langle C\rangle =\PP^5$.
\end{lemma}

\begin{proof}
Let $q$ be a generic conic. Assume that there exists a
reduced quartic $C_4^0$ 
passing through two distinct points $x,y$ of $q$, or
which is tangent to $q$ at one point $x=y$. 
We have
$l=\length (C_4^0\cap q)=2$, for if $l\geq 3$, then
$\deg \Psi_q(C_4^0)=2\deg C_4^0-3l=8-3l<0$, which is absurd.
In fact, the only irreducible curves $C\subset X$ whose degree with respect to
the linear system defining $\Psi_q$ is negative
are components of the
reducible members of the family $\CCC_3^0[2]_q$ contracted by
$\Psi_q$, so $\deg C\leq 2$.

The birational map
$\Psi_q$ transforms $C_4^0$ into
a conic meeting $\Gamma_q$ at 4 points, $u_1,u_2,u_3,u_4$. These points
span a plane $\PP^2$. As in the proof of Proposition \ref{lines},
consider $\Gamma_q$ as the projection of the canonical
curve $\Gamma$ from the line $\bar{uv}\subset\PP^6$, where $\lambda (q)=u+v$.
Then $\langle u_1,u_2,u_3,u_4,u,v\rangle=\PP^4$ and, by the
geometric Riemann--Roch Theorem, $\sum u_i+u+v\in G_6^1=G_6^1 (\Gamma )$,
where $G^r_d$ denotes the union of all
linear series $g_d^r$ on $\Gamma$; we keep the notation $W^r_d$
for the Brill--Noether locus of classes of such divisors in $\Pic^d(\Gamma )$.

Assume that $\Gamma$ (or equivalently, $X$) is generic. By \cite{ACGH},
$G_6^1$ is a $\PP^1$-bundle over
$W_6^1$, both $G_6^1$ and $W_6^1$ are nonsingular, irreducible,
and $\dim W_6^1=3$.

Thus we have constructed a map $\CCC_4^0[2]_q\lra G_q$,
where $G_q\subset G_6^1$ is the subset of divisors $D$
with $D-u-v$ effective. It is obvious that $G_q$ is 2-dimensional.
In fact, for generic $k\leq 4$ points $z_1,\ldots ,z_k\in\Gamma$,
the dimension of $G_{z_1,\ldots ,z_k}=\{ D\in G_6^1 
\mid D-\sum z_i\ \mbox{is effective}\}$ is equal to $4-k$.

It is easy to construct the inverse map: take
a divisor $D\in G_q$ and let $D-u-v=u_1+u_2+u_3+u_4$. Then,
after projecting to $\PP^4$ from $\bar{uv}$, we have
$\langle u_1,u_2,u_3,u_4\rangle_{\PP^4} =\PP^2$.
As $Q^3$ does not contain planes, $\PP^2\cap Q^3$
is a conic, say $C_2$, and $C_4^0:=\Psi_q^* (C_2)\in\CCC_4^0[2]_q$.
The scheme-theoretic intersection $C_4^0\cap q$ is either
two distinct points, or one point with multiplicity 2.

We have seen that $\CCC_4^0[2]_q$ is nonempty,
2-dimensional and birational to $G_q$.
Take another generic conic $q'$, and let $\lambda (q')=u'+v'$.
Then $G_q\cap G_{q'}=G_{u,v,u',v'}$ is finite. Hence
the union of $\CCC_4^0[2]_q$ when $q$ runs over an appropriate
open subset $U\subset\FFF (X)$ is 4-dimensional. 
This implies that the generic quartic from this union is
irreducible, for the family of reducible quartics
in $X$ is 3-dimensional. For the pairs of intersecting conics,
this follows from the fact that for a generic $x\in X$,
there are only finitely many (namely, 24) conics passing through
$x$, see Section 1. For the pairs of type a cubic plus a line, use Lemma \ref{cubics}.

We have seen that the family of good sextics of the form
$C_4^0+q$ is birational to $I^{(2)}$, where
$I^{(k)}=\{ (F,D)\in \Gamma^{(k)}\times G_6^1\mid
D-F\ \mbox{is effective}\}$.
The irreducibility of $I^{(2)}$ is proved in the same way
as in Lemma \ref{irred-quartics}. Denote by $q_k$
the natural projection to $G_6^1$ and restrict to
a generic pencil $\PP^1=g_6^1\subset G_6^1$. The
6-sheeted covering $q_1^{-1}(\PP^1)\lra\PP^1$ has only
simple ramifications, hence its monodromy is the whole
of $S_6$ and all the $I^{(k)}$ for $k=1,\ldots ,6$ are irreducible.

The fact that $C_4^0\cap q$ is generically a pair of distinct points
follows from the degeneration of $C_4^0$ to a curve of the form
$C_3^0+\ell$, where $\length (C_3^0\cap q)=2$, that is
$C_3^0\in\CCC_3^0[2]_q$ in the notation of Section \ref{sect-AJ}.
But the family $\CCC_3^0[2]_q$ is well understood: all its members
are smooth rational curves contracted by $\sigma_Q$,
except for 14 reducible members of the form $q_i+\ell_i$,
where $\ell_i$ are the flopping lines of $\Psi_q$, and
$q_i,\ell_i$ are unisecant to $q$. Hence the generic
$C_3^0\in\CCC_3^0[2]_q$ meets $q$ at two distinct points,
and the same is true for a generic $C_4^0\in\CCC_4^0[2]_q$.

Now, let $C$ be any good sextic in $X$.
Assume that the linear span of $C$ is strictly smaller than $\PP^5$.
Let, for example, $\langle C\rangle  = {\PP}^4$.
Then the projection from a general secant line
$<x,y>$, $x,y \in C$ sends $C$ to a quartic curve
$\overline{C} \subset {\PP}^2$ with at least two double points giving rise
to two 4-secant planes  to $C$ passing through
$<x,y>$. By Lemma \ref{4-secant}, these planes meet
$X$ along two conics passing through $x,y$, which contradicts
Lemma \ref{quartics}, (iii).
\end{proof}

\begin{proposition}\label{smoothings}
There is a distinguished
$6$-dimensional irreducible component $\CCC_6^{1*}(X)$
of the family
of elliptic sextics in $X$ satisfying the following
properties:

(i) The closure $\bar{\CCC}_6^{1*}(X)$ of $\CCC_6^{1*}(X)$ in $\Hilb_X$
contains the $4$-dimensional family $\CCC_{4,2}^1$ of reducible good
sextics of the form $C_4^0+q$
introduced in Lemma {\em\ref{good-sextics}}. 

(ii) A generic good sextic of the form $C_4^0+q$
is a smooth point of $\Hilb_X$.

(iii) A generic good sextic of the form $C_4^0+q$
can be partially smoothed to an irreducible rational
curve with only one node, and such partial smoothings fill
a five-dimensional subfamily of $\bar{\CCC}_6^{1*}(X)$.
\end{proposition}

\begin{proof}
For $C=C_1\cup C_2$ with $C_1=q$, $C_2=C_4^0$,
we have the following exact sequences \cite{HH}:
$$
0\to \NNN_{C/X}\to\bigoplus_{i=1}^2 N_{C/X}|_{C_i}\xrightarrow{\alpha}
\NNN_{C/X}|_Z\to 0,\ \ \ \
\length (\NNN_{C/X}|_Z)=4,
$$
$$
0\to \NNN_{C_i/X}\to \NNN_{C/X}|_{C_i}\xrightarrow{\varepsilon_i}T^1_Z\to 0,\ \ \
i=1,2,
$$
$$
0\to \NNN_{C/X}|_{C_i}(-Z)\to \NNN_{C/X}\xrightarrow{R_i}\NNN_{C/X}|_{C_{2-i}}\to 0,\
\ \ i=1,2,
$$
where $Z=C_1\cap C_2$, $\alpha$ is the difference map
$(s_1,s_2)\mapsto (s_2-s_1)|_Z$ and $T^1_Z$ is the Schlesinger
sheaf of infinetisimal deformations of singularities of $C$.
For generic $C$, $Z$ is a pair of distinct points and
$T_Z^1$ is a sky-scraper sheaf with 1-dimensional fibers
at points of $Z$, so that
$\length (T^1_Z)=2$.

A sufficient condition for the smoothness of
$\Hilb_X$ at $C$ is $h^1(\NNN_{C/X})=0$. If it is verified, then
the smoothability of $C$ is equivalent to the following condition:
the image of the composition
\begin{equation}\label{compo}
H^0(\NNN_{C/X})\xrightarrow{H^0R_i}H^0(\NNN_{C/X}|_{C_{2-i}})
\xrightarrow{H^0\varepsilon_{2-i}}H^0(T^1_Z)
\end{equation}
generates
the sheaf $T^1_Z$ for at least 
one value of $i$. The property (iii) is equivalent to
saying that one can
smooth by a small analytic deformation
only one node in a general curve of type $C_4^0+q$.
A sufficient condition which assures the existence
of such a partial smoothing is the surjectivity of
the map (\ref{compo}) for at least 
one value of $i$.

The three conditions are obviously verified if 
$\NNN_{C/X}|_{C_i}\simeq\mathcal{O}_{{\PP}^1}(a)\oplus\mathcal{O}_{{\PP}^1}(b)$
with $a>0,b>0$ for one value of $i$ and $a\geq 0,b\geq 0$ for the other.
The second exact sequence, Proposition \ref{conics}, (ii), and
Lemma \ref{quartics-normal}
imply that
$\NNN_{C/X}|_{C_1}\simeq\mathcal{O}_{{\PP}^1}(1)\oplus\mathcal{O}_{{\PP}^1}(1)$
or $\mathcal{O}_{{\PP}^1}\oplus\mathcal{O}_{{\PP}^1}(2)$
and
$\NNN_{C/X}|_{C_2}\simeq\mathcal{O}_{{\PP}^1}(2)\oplus\mathcal{O}_{{\PP}^1}(2)$ or
$\mathcal{O}_{{\PP}^1}(1)\oplus\mathcal{O}_{{\PP}^1}(3)$. This proves the
proposition. 
\end{proof}

%

\begin{corollary}\label{sextics-irred}
The family of elliptic sextics on a generic $X=X_{12}$
is irreducible: $\CCC_6^{1}(X)=\CCC_6^{1*}(X)$.
\end{corollary}

\begin{proof}
The proof is similar to that of Lemma \ref{irred-quartics}. A result of
\cite{P-2} is used, which states that the family of elliptic curves
$\CCC_d^1(\Sigma)$ of given degree $d\geq 4$ on the spinor tenfold
is irreducible. 
\end{proof}

\begin{lemma}\label{sextics-span}
For generic $C\in\CCC_6^{1}(X)$, \ $\langle C\rangle\cap X=C$
scheme-theoretically.
\end{lemma}

\begin{proof}
Notice that this is definitely false for some special
$C$, for there are elliptic sextics in $X$ having a secant
line contained in $X$. But one can show
that a generic elliptic sextic from 
${\CCC}_6^{1*}(X)$ has no secant lines. Indeed, if we assume the
contrary, then the generic quasi-elliptic sextic of the form $C_4^0+q$
has also a secant line, say $\ell$. This line is not a secant to $C_4^0$,
because by Lemma \ref{quartics-span}, $\langle C_4^0\rangle\cap X=C_4^0$
for a generic quartic $C_4^0$.
Hence $\ell$ is one of the 14 lines meeting $q$,
which are the flopping curves of $\Psi_q$.
Degenerate now $C_4^0$ to a curve of the form $C_3^0+\ell'$,
where $C_3^0\in\CCC_3^0[2]_q$ and $\ell'$ is a unisecant to $C_3^0$.
Then $\ell'$ is movable, hence generically different from $\ell$,
and both $C_3^0$ and $\ell'$ meet $\ell$. This is absurd, for
the generic member of $\CCC_3^0[2]_q$ is an exceptional curve
of $\sigma_Q$ which does not meet any of the flopping curves.

So, assume that $C$ has no secant lines and there is a
point $p\in \PP^5\cap X\smallsetminus C$. The 3-secant planes
$\PP^2$ of $C$ sweep over all the projective space
$\langle C\rangle =\PP^5$, so there is a 3-secant $\PP^2$ to $C$
passing through $p$. By Lemma \ref{4-secant}, there is a conic
$q$ in $X$ passing through the 4 points of $C\cap \PP^2$,
so $X$ contains the octic $C+q$ of arithmetic genus 3.
Except for $C$, $q$, there are no 
other curves in $\langle C\rangle\cap X$, for otherwise the residual
curves to $\langle C\rangle\cap X$ in the linear sections $\PP^6\cap X$
through $\langle C\rangle \cap X$ will form a rational net of
cubics, conics or lines in $X$, which is absurd. But in the
case when the 1-dimensional locus of $\langle C\rangle \cap X$
is $C\cup q$, we also obtain a contradiction:
the residual quartic curve $D$ in a generic $\PP^6$-section
of $X$ through $C+q$ satisfies $\length (D\cap (C+q))=5$.
As $D$ is reduced, it is a rational normal quartic,
hence $\langle D\cap (C+q)\rangle = \langle D\rangle =\PP^4$,
which is absurd, as $\langle C\rangle \cap \langle D\rangle=\PP^3$.

The above argument works as well if $p$ is an embedded
component of $\langle C\rangle\cap X$
whose tangent space is not contained in the tangent
space to the secant variety of $C$. Hence $C$ is a scheme-theoretic
intersection $\langle C\rangle\cap X$ for generic $C$.
\end{proof}

\section{The Abel--Jacobi map on elliptic sextics}

Let $X=\XT$ be a generic linear section $\sten\cap\PP^8$.
Exactly as in \cite{IM-3} in the case of quasi-elliptic
quintics, we can associate to any quasi-elliptic sextic $C\subset X$
a rank-2 vector bundle $\EEE =\EEE_C$ on $X$ with Chern
classes $c_1(\EEE )=H$ and $c_2(\EEE )=6[\ell ]$,
where $H$ is the class of a hyperplane section and $[\ell ]$
the class of a line.
It is obtained as the middle term of the following
nontrivial extension of $\OOO_X$-modules:
\begin{equation}\label{serre}
0\lra \OOO_X\lra \EEE \lra \III_C(1) \lra 0\; ,
\end{equation}
where $\III_C=\III_{C/X}$ is the ideal sheaf of $C$ in $X$.
One can easily verify (see \cite{MT-1} for a similar argument)
that, up to isomorphism, there is a unique
nontrivial extension (\ref{serre}), thus $C$ determines
the isomorphism class of $\EEE$. This way of constructing
vector bundles is called Serre's construction. The vector
bundle $\EEE$ has a section $s$ whose scheme of zeros
is exactly $C$. Conversely, for any section $s\in H^0(X,\EEE )$
such that its scheme of zeros $C_s=(s)_0$ is of codimension 2,
the vector bundle obtained by Serre's construction from
$C_s$ is isomorphic to $\EEE$. The normal sheaf $\NNN_{C_s/X}$
is naturally isomorphic to $\EEE|_{C_s}$.
As $\det\EEE\simeq\OOO_X(1)$, we
have $\EEE\simeq\EEE^\dual (1)$.

Let us denote by $M_X(2;m,n)$ the moduli space of stable
rank-2 vector bundles with fixed Chern classes
$c_i\in H^{2i}(X,\ZZ )$: $c_1=mH$ and $c_2=n[\ell ]$.
Recall also some notation from Section 5:
$\CCC_{4,2}^1(X)$,  the $4$-dimensional family of reducible good
sextics of the form $C_4^0+q$
introduced in Lemma {\ref{good-sextics}}, and $\CCC_6^{1}(X)$,
the
$6$-dimensional irreducible family of elliptic sextics in $X$.

\begin{lemma}\label{glob-gen}
For generic $C\in\CCC_6^{1}(X)$, the associated vector bundle $\EEE_C$
is generated by global sections.
\end{lemma}

\begin{proof}
By (\ref{serre}), it suffices to verify that $\III_C(1)$ is generated
by global sections, or equivalently, that
$\PP^5\cap X=C$ scheme-theoretically, where $\PP^5=\langle C\rangle $.
This follows from Lemma \ref{sextics-span}.
\end{proof}

The following proposition is proved in the same way as
similar statements for the (quasi-)elliptic quintics and associated
vector bundles in Section 3 of \cite{IM-3}.

\begin{proposition}\label{bundles}
For any good sextic $C\subset X$, the associated vector
bundle $\EEE$ possesses the following properties:

(i) $h^0(\EEE )=4$, $h^i(\EEE (-1))=0\ \forall\ i\in\ZZ$, and
$h^i(\EEE (k))=0\ \forall\ i>0$, $k\geq 0$.

(ii) $\EEE$ is stable and the local dimension of the
moduli space of stable vector bundles at $[\EEE ]$ is at
least $3$.

(iii) The scheme of zeros $(s)_0$ of any nonzero section
$s\in H^0(X,\EEE )$ is a quasi-elliptic sextic with
linear span $\PP^5$. 

(iv) If $s,s'$ are two nonproportional sections of $\EEE$,
then $(s)_0\neq (s')_0$. This means that $(s)_0$ and $(s')_0$
are different subschemes of $X$.

(v) The following three conditions are equivalent:\\
\hspace*{2 em}\newlength{\rrrrr}\rrrrr=\textwidth\addtolength{\rrrrr}{-2.2 em}
\begin{minipage}{\rrrrr}\par
(a)~for some (and hence for any) nonzero section $s'\in H^0(X,\EEE )$,
the Hilbert scheme of curves $\Hilb_X$ is nonsingular
and $6$-dimensional at $[C']$, where $C'=(s')_0$ is the zero locus
of $s'$;\par
(b)~the moduli space of stable rank-$2$
vector bundles $M_X(2;1,6)$ is nonsingular and $3$-dimensional at $\EEE$;
\par (c)~for some (and hence for any) nonzero section $s'\in H^0(X,\EEE )$,
$h^1(\NNN_{C'/X})=0$, where $C'=(s')_0$.
\end{minipage}
If, moreover, the zero loci $(s)_0$ for $s\in H^0(X,\EEE )$
have no base points,
then (a), (b), (c) are equivalent to:\\
\hspace*{2 em}
\begin{minipage}{\rrrrr}\par
(d) for some (and hence for any) nonzero section $s'\in H^0(X,\EEE )$,
$\NNN_{C/X}$ is a nontrivial extension of $\OOO_C$ by $\OOO_C(1)$,
that is, there is an exact triple
$$0\lra\OOO_C\lra\NNN_{C/X}\lra\OOO_C(1)\lra 0$$
and $\NNN_{C/X}\not\simeq \OOO_C\oplus\OOO_C(1)$.
\end{minipage}
\end{proposition}

The Serre construction can be relativized to provide a rational map
$\bar{\CCC}_6^{1}(X)\lrdash M_X=M_X(2;1,6)$, which we will call
the Serre map. 
Let $M^0_X$ be the image of the smooth
locus of ${\CCC}_6^{1}(X)$ in $M_X$
and $\CCC_X$ its inverse image in $\bar{\CCC}_6^{1}(X)$.
Propositions \ref{smoothings},
\ref{bundles} and Lemma \ref{glob-gen} imply the following corollary:

\begin{corollary}\label{serre-map}
(i)
$\CCC_X$, resp. $M^0_X$ is an open subset in the smooth
locus of $\bar{\CCC}_6^{1}(X)$, resp. $M_X$; $\dim\CCC_X=6$
$\dim M^0_X=3$, and the
Serre map ${\EuScript S}:\CCC_X\lra M^0_X$ is a
locally trivial \mbox{$\PP^3$-bundle}. 

(ii) $\CCC_X$ 
contains a $4$-dimensional family $\CCC_{4,2}^1\cap\CCC_X$ of reducible good
sextics of the form $C_4^0+q$.

(iii) The fiber ${\EuScript S}^{-1}([\EEE])\simeq
\PP H^0(X, \EEE)$ is identified with the family of zero loci $(s)_0$
of the sections $s\in H^0(X, \EEE)$ and consists of quasi-elliptic
sextics $C$ satisfying the condition $h^1(\NNN_{C/X})=0$.
\end{corollary}


Let now $\Gamma =\sten\cap\PP^6$ be the dual curve of genus 7
associated to $X$. The Brill--Noether
locus $W_6^1$ of $\Gamma$ is identified with the singular locus of
the canonical theta divisor $\Theta\subset \Pic^6(\Gamma )$ (see
\cite{GH}, Riemann--Kempf Theorem, Section 2.7).
Denote by
$\alp :\CCC_X\lra \Pic^{30}(\Gamma )$ the restriction of the canonical
Abel--Jacobi map to $\CCC_X$, $[C]\mapsto AJ(C)$ (see
Definition \ref{def-aj-can}).

\begin{theorem}\label{singtheta}
The Abel--Jacobi map $\alp :\CCC_X\lra \Pic^{30}(\Gamma )$
factors through the Serre map ${\EuScript S}$, that is
there exists a morphism \mbox{$\beta :M^0_X\lra \Pic^{30}(\Gamma )$}
such that $\alp =\beta\circ{\EuScript S}$. The map
$\beta$ is a birational isomorphism of $M^0_X$ onto
the singular locus of the divisor $3K -\Theta
\subset\Pic^{30}(\Gamma )$, where $K = K_\Gamma $ is the canonical
class of $\Gamma$.
\end{theorem}

\begin{proof}
The fibers of ${\EuScript S}$ are projective spaces, so they are
contracted to points by the Abel--Jacobi map. Thus $\beta$
exists as a set-theoretic map. The fact that it is a morphism
can be proved along the lines of the proof of Theorem 5.6 in \cite{MT-1}.

Consider the restriction of $\alp$ to the reducible sextics from
$\CCC_{4,2}^1$.
In the proof of Lemma \ref{good-sextics}, we described a birational isomorphism
of $\CCC_{4,2}^1$ with $I^{(2)}=\{ (F,D)\in \Gamma^{(2)}\times G_6^1\mid
D-F\ \mbox{is effective}\}$, where $G_6^1$ is the union of all the
linear series $g_6^1$ in $\Pic^{6}(\Gamma )$. Let $C_4^0+q$ be a generic
curve from $\CCC_{4,2}^1$, represented by a point $(F,D)\in I^{(2)}$.
By Lemma \ref{Tq} and Proposition \ref{aj-can},
$AJ(C_4^0)=4K-4\lambda (q)-[Z^q_{C_4^0}]-2d_8^q$. By construction,
$Z^q_{C_4^0}=D-F$, $\lambda (q)=[F]$. By Corollary \ref{aj-conics} and Lemma
\ref{d-eight}, $AJ(q)=K-\lambda (q)$ and $d_8^q=K-2\lambda (q)$.
This implies that $\alp (C_4^0+q)=AJ(C_4^0)+AJ(q)=3K-[D]$. As $D$ runs over
$G_6^1$, the classes $3K-[D]$ fill the divisor $3K-W_6^1=3K-\Sing\Theta$.
The image of $\alp$ coincides with that of $\beta$, and hence is
at most 3-dimensional, for $\dim M_X^0=3$. As $\dim W_6^1=3$,
$\alp (\CCC_{4,2}^1\cap\CCC_X)$ is dense in $M_X^0$ and
$\beta$ is quasifinite.

It remains to prove that $\beta$ is birational onto its image,
or equivalently, that the generic fiber of $\alp$ is one copy
of $\PP^3$. As the 4-dimensional family $\CCC_{4,2}^1\cap\CCC_X$
dominates $M_X^0$, the fibers $\PP^3$ of ${\EuScript S}$ contain
generically a 1-dimensional family of curves from $\CCC_{4,2}^1$.
So, if there were several fibers of ${\EuScript S}$ in one fiber
of $\alp$, then the generic fiber of the restriction
$\CCC_{4,2}^1\cap\CCC_X\lra 3K-W_6^1$ of $\alp$ would be a disjoint union of
several curves. But we have seen in the proof of Lemma \ref{good-sextics}
that this fiber is an irreducible 15-sheeted covering of $\PP^1$,
so the generic fiber of $\alp$ is connected.
\end{proof}

\

\section{Irreducibility of $M_X(2;1,6)$}

Let $X=\XT$ be a Fano 3-dimensional linear section 
of the spinor tenfold and $M_X=M_X(2;1,6)$. 
We will prove that $M_X$ is irreducible for generic $X$.
This will follow from the irreducibility of the family of elliptic sextics
on a generic $X$ as soon as we have proved
that a generic  $\EEE$ in any component of $M_X$
is obtained by Serre's construction from an elliptic sextic.

\begin{lemma}\label{E-from-Z}
Let
$\EEE\in M_X$, $S$ a generic hyperplane section of $X$,
$E=\EEE|_S$ the restriction of $\EEE$ to $S$.
Then the following assertions hold:

(i) $\chi (\EEE)=4$, $h^3(\EEE )=0$.

(ii) $E$ is stable and the scheme of zeros $Z_s=(s)_0$ of 
any nonzero section $s$
of $E$ is $0$-dimensional and of length $6$.
$E$ can be obtained by Serre's construction on $S$ from
a $0$-dimensional subscheme $Z\subset S$ of length $6$:
\begin{equation}\label{SerreZ}
0\lra \OOO_S \lra E \lra \III_Z(1) \lra 0.
\end{equation}
For such a $Z$, $\dim\langle Z\rangle=4$
and $\langle Z\rangle=\langle Z'\rangle$ for any $Z'\subset Z$
of length $5$.

(iii) $E$ is generated by global sections at the generic point
of $S$. 
\end{lemma}

\begin{proof}
(i) We have $\chi (\EEE)=4$ by Riemann--Roch,
and $h^3(\EEE)=h^0(\EEE (-2))=0$ by stability.

(ii) $E=\EEE|_S$ is slope-semistable by Theorem 3.1 of \cite{Ma}.
The semistability implies the stability because $\Pic S=\ZZ H$ and
$\det E=\OOO(H)$ is odd. Hence
$h^2(E)=h^0(E (-1))=0$ and $\chi (E)=4$ implies
$h^0(E)\geq 4$. The zero locus $(s)_0$ of any non-zero
section of $E$ is finite, for otherwise it would be a curve
from the linear system $|kH|$ and then $h^0(E(-k))\neq 0$,
which is absurd. Hence it is a subscheme $Z$ of length equal to $c_2(E)=6$,
and there is an exact triple (\ref{SerreZ}) with the inclusion
$\OOO_S \lra E$ defined by~$s$. We have $h^1(\III_Z(1))=5-m$,
where $m=\dim\langle Z\rangle $. By Serre duality,
$\dim \Ext^i(\III_Z(1),\OOO_S)=h^{2-i}(\III_Z(1))$, hence
the triple (\ref{SerreZ}) can be nonsplit only if $m\leq 4$.
The values $m\leq 2$ are impossible by Lemma \ref{4-secant}.
Hence $m=3$ or 4.

Assume that $m=3$. By Lemma \ref{4-secant}, 
for any subscheme $Z'\subset Z$ of length
5, we have $\dim\langle Z'\rangle= \dim\langle Z\rangle =3$. 
Hence $h^1(\III_{Z'}(1))=1$, and
there is a unique nontrivial extension 
$$0\lra \OOO_S \lra E' \lra \III_{Z'}(1) \lra 0.$$
Again by Lemma \ref{4-secant}, for any $Z''\subset Z'$ of length
4, $\langle Z''\rangle=\langle Z'\rangle$, which implies the
local freeness of $E'$ (see, for example, \cite{Tyu}, Lemma 1.2).
Thus the Serre construction
applied to $Z'$ provides a rank-2 vector bundle $E'$
with $c_1(E')=[H]$, $c_2(E')=5$. It is easy to see that
$E'$ is stable. Indeed, if we assume that it
is unstable, then any destabilizing subsheaf
should be of the form $\III_{W}(k)$, where $k>0$ and
$W$ is a 0-dimensional subscheme of $S$. If we replace
$\III_{W}(k)$ by its saturate $\III_{W}(k)^{\dual\dual}=\OOO_S (k)$,
we get an inclusion
$\OOO_S(k)\into E'$, which is absurd, since $h^0(\OOO_S(k))\geq
h^0(\OOO_S(1))=8>h^0(E')=5$. By Corollary 5.8 of \cite{IM-3},
$E'$ is generated by global sections. From Serre's exact
triple for $E'$, we conclude that $\III_{Z'}(1)$ is generated
by global sections. Hence $\langle Z'\rangle\cap S=Z'$
scheme-theoretically, which contradicts the equality
$\langle Z'\rangle= \langle Z\rangle$. Thus we
have proved that $m=4$, that is, $\langle Z\rangle\simeq \PP^4$.

Suppose now that there is a subscheme $Z'\subset Z$ of length 5
with $\langle Z'\rangle\subsetneqq \langle Z\rangle$. Then
$\dim \langle Z'\rangle=3$. The exact triple
$
0\rar\III_Z(1)\xrightarrow{\iota}\III_{Z'}(1)\rar\CC_p\rar 0,
$
where $\{p\}$ is the support of $\III_{Z'}/\III_Z$, and the
local-to-global spectral sequence provide the
commutative diagram
$$
\xymatrix{
\Ext^1(\III_Z(1),\OOO_S)\ar@{^(->}[r] &
    H^0(\EXT^1(\III_Z(1),\OOO_S)) \\
\Ext^1(\III_{Z'}(1),\OOO_S)\ar@{^(->}[r]\ar[u]^{\simvert}_{\iota^*}&
H^0(\EXT^1(\III_{Z'}(1),\OOO_S))\ar@{^(->}[u]
}
$$ 
>From its right column we see that the extension class of
(\ref{SerreZ}) does not generate the stalk of
$\EXT^1(\III_Z(1),\OOO_S)$, which contradicts the local freeness
of $E$ at $p$ by Serre's Lemma (see e. g. Lemma  5.1.2 in \cite{OSS}).
Hence $\langle Z'\rangle= \langle Z\rangle$ and we are done.

(iii) Let $s_1,s_2$ be two non-proportional sections of $E$.
If they do not generate $E$ at any point of $S$,
then there is a rank-1 subsheaf of $E$ with at least 2
linearly independent sections, which contradicts the stability.
\end{proof}

\begin{lemma}\label{E-cohomology}
In the assumptions of Lemma \ref{E-from-Z}, the following statements
hold:

(i) $h^1(E(k))=0$ for all $k\in\ZZ$ and $\chi (E(k))=h^0(E(k))=h^2(E(-k-1))=12k(k+1)+4$
for $k\geq 0$.

(ii) $h^i(\EEE (k))=0$ for all $k\in\ZZ$, $i=1,2$;
$\chi (\EEE (k)) =h^0(\EEE (k))=h^3(\EEE (-k-2))=4(k+1)^3$ for $k\geq -1$.
\end{lemma}

\begin{proof}
This is standard; use the exact triples
$$0\lra \OOO_S (k)\lra E (k) \lra \III_Z(k+1) \lra 0,$$
$$0\lra \III_Z(k) \lra \OOO_S (k)\lra \OOO_Z (k)\lra 0,$$
$$0\lra \EEE (k-1)\lra \EEE (k)\lra E (k)\lra 0,$$
the Serre duality and the Kodaira vanishing $h^1(\EEE (k))=0$
for \mbox{$k\ll 0$}\sloppy .
\end{proof}

\begin{corollary}\label{EEE-cohomology}
Let
$\EEE\in M_X$, $S$ any nonsingular hyperplane section of $X$,
$E=\EEE|_S$ the restriction of $\EEE$ to $S$. Then the
restriction map $H^0(\EEE)\lra H^0(E)$ is an isomorphism
and the assertion (i) of Lemma \ref{E-cohomology} holds
for the cohomology $h^i(E (k))$.

If in addition $\Pic S= \ZZ H$, then $E$ is stable and any two
nonproportional sections of $E$ define distinct $0$-dimensional
length-$6$ subschemes of $S$.
\end{corollary}

\begin{proof} The assertions on the restriction map
and on $h^i(E (k))$ are obvious.
If $\Pic S= \ZZ H$, $\rk E=2$ and $c_1(E)=H$, then 
the stability of $E$ is equivalent to
$h^0(E (-1))=0$. Hence $E$ is stable. This imples that
$\Hom (E,E)=H^0(E\otimes E(-1))=\CC$. From the exact triple
(\ref{SerreZ}) tensored by $E(-1)$, we deduce that
$H^0(E\otimes\III_Z)\simeq H^0(E\otimes E(-1))=\CC$,
hence a section of $E$ having $Z$ as its zero locus
is unique up to proportionality.
\end{proof}

\begin{proposition}\label{E-sextic}
Assume $X$ generic, and
let $\EEE\in M_X$. Then $\EEE$ is generated by global sections
at the generic point of $X$ and
can be obtained by Serre's construction from a quasi-elliptic sextic
lying in the closure of the family of elliptic sextics
in the Hilbert scheme of $X$.
\end{proposition}

\begin{proof}
Let $C_s$ denote the curve $(s)_0$ of zeros of
a nonzero section $s\in H^0(\EEE)$. It is a
l.~c.~i. sextic curve with trivial canonical sheaf for any
$s\neq 0$. Moreover, it is connected,
that is $h^0(\OOO_{C_s})=1$, and $\NNN_{C_s}\simeq\EEE|_{C_s}$,
so $\chi (\NNN_{C_s})=6$. This implies that the dimension
of the Hilbert scheme of curves in $X$ at the point
$\{C_s\}$ representing $C_s$ is at least 6. Moreover,
the properties of being a l.~c.~i. curve and to have
trivial canonical sheaf are open, so
any small deformation of a l.~c.~i. curve with trivial canonical sheaf
is of the same type.
We will use this
observation to show that $C_s$
is in the closure of the family of smooth elliptic sextics in $X$.

The outline of the proof is the following. First,
we decompose $C_s$
into the sum of the fixed part $F$ and the movable part $M_s$. Second, we
show that $\deg M_s\geq 4$. Finally, in assuming $s$ generic, we
examine the possible types of
decomposition of $M_s$ and $F$ in irreducible components to show that
$F+M_s$ deforms to a smooth sextic curve.
 
By Lemma
\ref{E-cohomology} and Corollary \ref{EEE-cohomology},
the curves $C_s$ form a family with base $\PP^3$,
and for two nonproportional sections $s,s'$ of $\EEE$,
the curves $C_s$, $C_{s'}$ are distinct as subschemes in $X$.
Let $F$ be the sum of the fixed components of this family, and $M_s$
the movable part, so that $C_s=F+M_s$ as an algebraic cycle.
By Bertini Theorem, both $F$ and the singular
loci of $M_s$ for generic $s$ (if nonempty) are contained in the base
locus BL$(\EEE)$ of $\EEE$, defined as the locus of points $x\in X$ in which
the stalk $\EEE_x$ is not generated by $H^0(\EEE)$. According to Lemma
\ref{E-from-Z}, (iii) and Corollary \ref{EEE-cohomology},
BL$(\EEE)$ is a proper closed subset of $X$,
so $M_s$ is reduced for generic $s$. Taking any 3 nonproportional sections
$s_1,s_2,s_3$ of $\EEE$ and a generic point $x\in X$, we can find
a nontrivial linear combination $s=\lambda_1s_1+\lambda_2s_2+\lambda_3s_3$
vanishing at $x$. Hence the family $\{M_s\}_{[s]\in\PP^3}$ is a covering family
of curves on $X$: there is at least one curve $M_s$ passing through
a generic point of $x$.

Let us show that the curves $M_s$
are different for nonproportional sections $s$. If $C_s$ has multiple components,
this does not follow directly from
the above observation that $C_s$, $C_{s'}$ are distinct whenever $s\not\sim s'$, 
for $C_s$, $C_{s'}$ may
differ, a priori, by the nilpotent structure along the multiple components whilst
the associated algebraic cycles $F+M_s$, $F+M_{s'}$ are the same. Thus, 
assuming that $s\not\sim s'$, we will verify that the supports
of $C_s$ and $C_{s'}$ are distinct.

By the stability of $\EEE$, the subsheaf 
$\OOO_X\cdot s+\OOO_X\cdot s'\subset\EEE$
cannot be of rank 1. Hence $s,s'$ are generically linearly independent and
the section $s\wedge s'\in H^0(\det \EEE)=H^0(\OOO_X(1))$ is nonzero.
As $\Pic X\simeq \ZZ$,
the zero locus $S=(s\wedge s')_0$ is a possibly singular, but reduced
and irreducible surface from the linear series of hyperplane sections of $X$.
Obviously $C_s\subset S$, $C_{s'}\subset S$.
The restrictions $\sigma=s|_S$, $\sigma'=s'|_S$ are 
sections of the rank-1 torsion-free sheaf
$\LLL=\OOO_S\cdot \sigma+\OOO_S\cdot \sigma'\subset
\EEE|_S$. They are nonproportional, 
for if $\lambda \sigma+\lambda' \sigma'=0$ for some
nonzero constants $\lambda,\lambda'\in\CC$, then
$\lambda s+\lambda' s'$ is a nonzero section of $\EEE$
which vanishes exactly on $S$,
and this is impossible by the stability of $\EEE$.

If we assume that $\Supp C_s=\Supp C_{s'}$,
 then all the nontrivial linear combinations
$\lambda s+\lambda' s'$ have the same zero set. This is absurd, for if
$x\in S$ is generic, then the fiber $\LLL(x):=\LLL\otimes \CC(x)$ is one-dimensional,
so there exists a nontrivial linear combination $\lambda s+\lambda' s'$
vanishing at $x$ and the curve $C_{\lambda s+\lambda' s'}$ passes through $x$.
This implies that $M_{\lambda s+\lambda' s'}$ is a movable curve,
and hence $M_s\neq M_{s'}$.

Since the family of lines is not covering for $X$
and since the one of conics contains no rational subvarieties, we
have $\deg M_s\geq 3$. Suppose that $\deg M_s= 3$.
Then we get a 3-dimensional family of cubic curves $\MMM=\{M_s\}$, bijectively
parameterized by $\PP^3=\PP H^0(\EEE)$. Let $s\in \PP H^0(\EEE)$ be generic. 
We have seen that then $M_s$ is reduced. Let us show that it is
also irreducible. Indeed, if $M_s$ is a line plus a conic, then
by projecting $\MMM$ to the families $\tau(X)$, $\FFF (X)$
of lines and conics in $X$,
we get a nonconstant rational map $\PP^3\dasharrow \tau(X)\times\FFF (X)$.
But $\tau(X)$ is a smooth curve of genus 43, and $\FFF(X)\simeq \Gamma ^{(2)}$
for the orthogonal genus-7 curve $\Gamma=\check{X}$,
so $\tau(X)$, $\FFF (X)$ do not contain
rational subvarieties. Further, if $M_s$ is a union of three lines,
we get a generically injective rational map $\PP^3\dasharrow \tau(X)^{(3)}$
which is also absurd, since
$\tau(X)^{(3)}$ is irreducible and nonrational.

Thus $M_s$
is a reduced and irreducible cubic curve in $X$ for generic $s$. As 
$X$ is an intersection of quadrics, the span of $M_s$ is $\PP^3$
and $M_s$ is nonsingular. We obtain a family of rational normal 
cubics in $X$, bijectively parameterized by an open set of $\PP^3$.
This contradicts Lemma \ref{cubics},
saying that $\CCC^0_3(X)$ is irreducible and birational to $\Gamma ^{(3)}$.

We have proved that $\deg M_s\geq 4$.  From now on we assume
$s\in H^0(\EEE)$ generic. We will treat several cases differing by the degree
of $M=M_s$ and the type of its decomposition in irreducible components. 

Case 1: $\deg M=6$, that is, $F=0$.
Then $C=M$ is a good sextic.

Subcase 1.1: $M$ is irreducible.
Either it is an elliptic sextic, and we are done, or
it is a rational sextic with one double point whose contribution
to the arithmetic genus is 1, that is a node or a cusp. An argument as
in the proof of Lemma \ref{irred-quartics} shows that when
$X$ is generic, then all the components of the family of rational sextics
in $X$ are 6-dimensional, and the singular rational sextics 
fill a codimension-1 locus. As the local deformation space of
$C$ in $X$ is at least 6-dimensional, we conclude that
$C$ deforms to a smooth elliptic sextic in $X$.

Subcase 1.2: at least one of the components of  $M$
is a line. By the same argument as we used for $\deg M=3$,
the number of line components is $\geq 4$. But then the remaining component
cannot be a conic, for then this conic should be fixed
and $\deg M=4$, which is absurd. So, $M$ has to be
a connected union of 6 lines. As any line
meets only finitely many lines in $X$, the dimension
of the family of connected unions of 6 lines in $X$ is $\leq 1$,
but the family of different $M$'s is 3-dimensional,
so $M$ is not of this type.

Subcase 1.3: $M$ has a conic component. Then, as above,
$M$ is a connected union of three smooth conics, $M=q_1\cup q_2\cup q_3$.
The sextuples $\sum\lambda (q_i)$ of points of $\Gamma$, where
$\lambda:\FFF(X)\isoto \Gamma^{(2)}$ was defined in the proof
of Proposition \ref{conics},
sweep out a unirational 3-dimensional subvariety of $\Gamma^{(6)}$.
This is impossible, for $\Gamma$ is a generic genus-7 curve and hence it has
no $g^3_6$ (and even $g^2_6$).

Subcase 1.4: $M$ has a cubic component. Then it is a union of
two rational normal cubics $C_1\cup C_2$ with $\length (C_1\cap C_2)=2$.
As $\dim \CCC^0_3(X)=3$ and
two generic rational normal cubics in $X$ are disjoint, we see that the family
of pairs of intersecting rational normal cubics is at most 5-dimensional.
Hence $M$ deforms to an irreducible sextic, and this reduces the problem
to Subcase 1.1, which we have already settled.

Case 2: $\deg M=5$, then $F=\ell$ is a line. Similarly to the above,
we can prove that $M$ is a smooth
rational quintic and $\length (\ell\cap M_s)=2$. An argument as
in the proof of Lemma \ref{irred-quartics} shows that the rational quintics
in a generic $X$ fill a 5-dimensional family, and those
meeting a line twice lie in codimension 1. Hence $\ell+M$
deforms to an irreducible good sextic, which brings us to Subcase 1.1.

Case 3: $\deg M=4$. We have two subcases. 

Subcase 3.1: $C=q+M$,
where $q$ is a reduced conic. Then the result follows by the same
argument as in Case 2. 

Subcase 3.2: $C=F+M$,
where $F$ is a Cohen--Macaulay double structure on a line $\ell$. As before,
we can prove that $M$ is a rational normal quartic such that $\length (F\cap M)=2$.
We have an exact triple 
$$
0\lra \OOO_F(-Z)\lra\OOO_C\lra\OOO_M\lra 0,
$$
where $Z$ is the
intersection scheme $F\cap M$.

The Cohen--Macaulay double structures on a smooth curve are completely
described, for example, in \cite{BanF}. They all are obtained by
Ferrand's construction as in the proof of Lemma \ref{double-lines};
one can think of $F$ as $\ell$ together with a cross section $\xi$
of the projectivized normal bundle $\PP(\NNN_{\ell/X})$ over $\ell$.
The multiplicity of the intersection $F\cap M$ can be interpreted via
the relative position of the proper transform $\tilde{M}$ of $M$ and $\xi$
on the blowup $\tilde{X}$ of $X$ with center $\ell$.
We have the following three possibilities for the intersection
$F\cap M$ of total multiplicity~2:
(1) $\ell$ intersects $M$ quasi-transeversely at one point $p$, and $\tilde{M}$
passes through $\xi_p$; (2) $\ell\cap M=\{p_1,p_2\}$,
$p_1\neq p_2$,  and $\tilde{M}$
does not pass through any one of the points $\xi_{p_1},\xi_{p_2}$; 
(3) $M$ is simply tangent to $\ell$ at $p$ and $\tilde{M}$
does not pass through $\xi_p$. The singular points
of $C$ are analytically equivalent to $(x^2,z)\cap (x,y)$
in the case (2) and $(x^2,z)\cap (x, z-y^2)$ in the case (3).
These singularities are not Gorenstein, so the only
possible case is (1). But in this case we have
$\omega_F\simeq\omega_C|_F(-Z)$ and $\omega_C\simeq\OOO_C$.
Restricting $\omega_F$ to $\ell$, we obtain a contradiction
as in the proof of Lemma \ref{double-lines}: on the one hand,
$\omega_F|_\ell=\OOO_F(-Z)|_\ell=\OOO_{\ell}(-1)$,
on the other hand $\omega_F|_\ell\simeq\omega_\ell\otimes\LLL^{-1}$, where
$\LLL\simeq\OOO_\ell(k)$ for some $k\geq 0$, which is impossible.
Thus the Subcase 3.2 does not occur. 
\end{proof}

\begin{corollary}\label{MX-irred}
If $X$ is generic, then $M_X$ is irreducible.
\end{corollary}

\begin{proof}
This is an immediate consequence of Proposition \ref{E-sextic}
and Corollary \ref{sextics-irred}.
\end{proof}

\section{Appendix: Maps from the symmetric square of a curve}

Here we prove the following assertion: 

\begin{proposition}\label{aut}
Let $\Gamma$ be a generic curve of genus $g\geq 4$
and $S= \Gamma^{(2)}$ the symmetric square of $\Gamma$. Then 
the following assertions hold:

(i) If $g\neq 4$, then for any nonrational irreducible
curve $C$, there are
no nonconstant rational maps $\phi :S\dasharrow C$.

(ii) Let 
$\phi :S\dasharrow S$ be a nonconstant rational map. Then
$\phi=\id_S$.
\end{proposition}

We fix for the sequel the notations $\Gamma$ and $S$ for a
{\em generic} curve of genus $g$ and its symmetric square respectively.
The proposition follows from a sequence of lemmas. Before stating
them, we need to describe the Mori cone and the ample cone of $S$.

Let $g\geq 2$.
Let $\pi:\Gamma \times\Gamma \lra S=\Gamma^{(2)}$ be the
quotient map and $\Delta\subset \Gamma \times\Gamma $ the diagonal.
The Neron--Severi group $NS (S)$ contains 3 natural classes:
the first one is $f$, the class of a fiber $\pi (\{ x\}\times\Gamma )$,
where $x\in\Gamma $
is a point, the second one is $\delta =\frac{1}{2}\pi (\Delta )$,
and the third one is $\Theta|_S$, the pullback of the theta-divisor
via the Abel-Jacobi map $S\lra J(X)$ defined up to translations.
There is one relation among them, $\delta =(g+1)f-\Theta|_S$,
and $NS (S)$ is freely generated by $f$ and $\delta$
(see \cite{ACGH}, Sect. 5 of Ch. VIII, and \cite{GH}, Sect. 5 of Ch. II).
We have also: 
$$
\delta^2=1-g,\ \delta f=f^2=1,\ K_S=-\delta +(2g-2)f, \ 
K_S^2=(2g-3)^2-g,\ $$
$$c_2(S)=(2g-3)(g-1),\ \chi (\OOO_S)=
\frac{(g-1)(g-2)}{2}.
$$

If $g\geq 3$, then $S$ contains no rational
curves and $K_S^2>0$, so $S$ is of general type, and moreover,
$K_S$ is ample.

Let $N(S)$ be the real vector plane $NS(S)\otimes\RR$,
$\NE (S)\subset N(S)$ the smallest closed cone containing
the classes of effective curves (the Mori cone of $S$),
and $\NA (S)$ the dual cone with respect to the intersection
product on $N(S)$; this is the smallest closed cone containing
the classes of ample curves. In our case, the cones are just
angles in the plane. It is obvious that one of the rays
bordering $\NE (S)$ is $\RR_+\delta$ and the other is
of the form $\RR_+(-\delta+kf)$ for some real $k$, $1<k< g+1$.
Similarly, $\NA (S)$ is bordered by the rays $\RR_+(\delta+(g-1)f)$
and $\RR_+(-\delta+lf)$ with $k\leq l=\frac{k+g-1}{k-1}<g+1$
The following theorem, proved in
\cite{Kou}, \cite{CiKou}, gives more precise estimates:

\begin{thm}[Kouvidakis, Ciliberto--Kouvidakis]
Assume that $\Gamma$ is a generic curve
of genus $g\geq 4$. Then 
$$\sqrt{g}\leq k\leq \sqrt{g}+1\leq l
\leq\frac{g}{\sqrt{g}-1}+1.$$

If $\sqrt{g}\in\ZZ$, then $k=l=\sqrt{g}+1$. If moreover $g\neq 4$,
then there are no classes of effective
curves in the ray $\RR_+(-\delta+(\sqrt{g}+1)f)$.
\end{thm}



\begin{lemma}
Let $g\geq 5$, and let $C$ be a nonsingular complete
curve. Then there are no
nonconstant morphisms $\phi :S\lra C$. If moreover
$C$ is nonrational, then every rational map $\phi :S\dasharrow C$
is regular, hence constant.
\end{lemma}

\begin{proof}
The fiber of such a morphism would provide a rational numerically
effective class $h=a\delta +
bf$ with $h^2=0$, which implies $\frac{b}{a}=-1\pm\sqrt{g}$.
Hence $\sqrt{g}\in\ZZ$. In the interior of $\NA (S)$, $h^2>0$,
hence $h$ is on the border and is proportional to
$h_0=-\delta+(\sqrt{g}+1)f$. This contradicts
the non-existence of effective curves in the ray $\RR_+h_0$.
\end{proof}

\begin{remark}
For $g=4$, $\Gamma$ has two $g_3^1$'s. A $g_3^1$ defines
the following curve on $S$:
$$
D=\{ x+y\in S\mid \exists z\in\Gamma : \ x+y+z\in g_3^1\} .
$$
The two $g_3^1$'s thus provide
two curves $D,D'$ in $S$ in the same numerical class $-\delta +3f$
such that $D^2=D^{\prime 2}=0$. Hence the border ray of $\NE (S)$ contains
effective curves, and to extend the previous lemma
to $g=4$, one has to show that $\dim |nD|=\dim |nD'|=0$ for all
$n>0$.
\end{remark}

\begin{lemma}
Let $g\geq 3$, and let $\phi :S\dasharrow S$ be a rational map of
degree $d>0$. Then $d=1$.
\end{lemma}

\begin{proof}
Since $C$ is not hyperelliptic, $\phi$ is regular.
By \cite{Beau-2}, Proposition 2, if a compact complex manifold
$X$ admits an endomorphism of degree $d>1$, then $\kappa (X)
<\dim X$. But $S$ is of general type, so it has no endomorphisms
of degree $>1$.
\end{proof}

\begin{lemma}
Let $g\geq 4$ and let $\phi :S\dasharrow S$ be a birational map.
Then $\phi =\id_S$.
\end{lemma}

\begin{proof}
As $S$ contains no rational curves, $\phi$ is biregular.
The induced automorphism $\phi^*$ of $N(S)$ is given by an
integer matrix in the basis $\delta,f$. It preserves the
intersection product and the cones $\NE (S)$, $\NA (S)$.
The canonical class $K_S$ is an eigenvector of $\phi^*$
with eigenvalue 1. Hence if $\phi^*$ preserves both
border rays of $\NE (S)$, it is the identity map.
As $\delta^2<0$, there is
only one effective curve in the numerical class $2\delta$,
the diagonal $\Delta '=\pi (\Delta )$, so $\Delta'$ is invariant
under $\phi$. But $\Delta ' \simeq \Gamma$ and $\Gamma$ has no
nontrivial automorphisms, for it is a generic curve of genus $g$.
Hence $\phi|_{\Delta '} =\id$. 

Now, any of the curves $F_x=
\pi (\{ x\}\times\Gamma )$, represented by the class $f$,
is tangent to $\Delta'$ at a single point $2x=\pi (x,x)$.
Hence its image $\phi (F_x)$ is also tangent to $\Delta'$ at
$2x$ and belongs to the same class $f$. Lifting it to
$\Gamma\times\Gamma$, one immediately verifies that
$\phi (F_x)=F_x$ and $\phi =\id$.

It remains to consider the second case, when $\phi^*$
permutes the border rays of $\NE (S)$. Then
$\phi^*$ is an orthogonal reflection with mirror
$\RR K_S$. We have
$$
\phi^*(v)=v-2\frac{(v,\alp )}{(\alp ,\alp)}\alp,\ \ \alp=-(2g-3)
\delta + (3g-3)f.
$$
This gives $\phi^*(\delta )=-\frac{4g-3}{4g-9}\delta
+\frac{12g-12}{4g-9}f$. The coefficient of $\delta$ is
fractional for all $g\geq 4$, which contradicts the condition that
$\phi^*$ is integer in the basis $\delta, f$. Hence
the second case is impossible.
\end{proof}

\begin{remark}
The previous lemma does not extend to $g=3$,
because in this case the formula $\phi : x+y\mapsto K_\Gamma-x-y$
defines an involution on $S$.
\end{remark}


\end{document}